\documentclass[a4paper]{article}

\usepackage[cp1250]{inputenc}
\usepackage[IL2]{fontenc}
\usepackage{a4wide}
\usepackage[english]{babel}
\usepackage{euscript}
\usepackage{amstext,amsbsy,amscd,amssymb}
\usepackage{amsmath,enumerate}
\usepackage{amsfonts}
\usepackage{graphics}
\usepackage{graphicx}
\usepackage{microtype}
\usepackage{color}
\include{diagxy}
\allowdisplaybreaks

\let\rarr=\rightarrow

\let\veps=\varepsilon
\let\mcal=\mathcal
\let\mfrak=\mathfrak
\let\eus=\EuScript

\def\N{\mathbb{N}}

\def\C{\mathbb{C}}

\def\End{\mathop {\rm End} \nolimits}

\def\Ext{\mathop {\rm Ext} \nolimits}

\def\ad{\mathop {\rm ad} \nolimits}

\def\id{\mathop {\rm id} \nolimits}

\def\Sing{\mathop {\rm Sing} \nolimits}
\def\ch{\mathop {\rm ch} \nolimits}
\def\Cas{\mathop {\rm Cas} \nolimits}

\newsymbol\squares 1003

\long\def\proof #1{\noindent \emph{Proof.}\ #1 \hfill $\squares$
\medskip}

\newcounter{num}[section]
\numberwithin{equation}{section}
\numberwithin{num}{section}

\long\def\definition #1 {\refstepcounter{num} \noindent {\bf
Definition \thenum.} #1

\medskip}

\long\def\theorem #1{\refstepcounter{num} \noindent {\bf Theorem
\thenum.} #1

\medskip}

\long\def\lemma #1{\refstepcounter{num}  \noindent {\bf Lemma
\thenum.} #1

\medskip}

\long\def\proposition #1{\refstepcounter{num}  \noindent {\bf
Proposition \thenum.} #1

\medskip}

\newenvironment{enum}{\begin{list}{}{\topsep=2pt \itemsep=0pt
\parsep=0pt}}{\end{list}}

\newcommand*\riso{%
  \xrightarrow[]{\raisebox{-0.25em}{\smash{\ensuremath{\sim}}}}%
}

\makeatletter
\newcommand*\if@single[3]{%
  \setbox0\hbox{${\mathaccent"0362{#1}}^H$}%
  \setbox2\hbox{${\mathaccent"0362{\kern0pt#1}}^H$}%
  \ifdim\ht0=\ht2 #3\else #2\fi
  }
\newcommand*\rel@kern[1]{\kern#1\dimexpr\macc@kerna}
\newcommand*\widebar[1]{\@ifnextchar^{{\wide@bar{#1}{0}}}{\wide@bar{#1}{1}}}
\newcommand*\wide@bar[2]{\if@single{#1}{\wide@bar@{#1}{#2}{1}}{\wide@bar@{#1}{#2}{2}}}
\newcommand*\wide@bar@[3]{%
  \begingroup
  \def\mathaccent##1##2{%
    \if#32 \let\macc@nucleus\first@char \fi
    \setbox\z@\hbox{$\macc@style{\macc@nucleus}_{}$}%
    \setbox\tw@\hbox{$\macc@style{\macc@nucleus}{}_{}$}%
    \dimen@\wd\tw@
    \advance\dimen@-\wd\z@
    \divide\dimen@ 3
    \@tempdima\wd\tw@
    \advance\@tempdima-\scriptspace
    \divide\@tempdima 10
    \advance\dimen@-\@tempdima
    \ifdim\dimen@>\z@ \dimen@0pt\fi
    \rel@kern{0.6}\kern-\dimen@
    \if#31
      \overline{\rel@kern{-0.6}\kern\dimen@\macc@nucleus\rel@kern{0.4}\kern\dimen@}%
      \advance\dimen@0.4\dimexpr\macc@kerna
      \let\final@kern#2%
      \ifdim\dimen@<\z@ \let\final@kern1\fi
      \if\final@kern1 \kern-\dimen@\fi
    \else
      \overline{\rel@kern{-0.6}\kern\dimen@#1}%
    \fi
  }%
  \macc@depth\@ne
  \let\math@bgroup\@empty \let\math@egroup\macc@set@skewchar
  \mathsurround\z@ \frozen@everymath{\mathgroup\macc@group\relax}%
  \macc@set@skewchar\relax
  \let\mathaccentV\macc@nested@a
  \if#31
    \macc@nested@a\relax111{#1}%
  \else
    \def\gobble@till@marker##1\endmarker{}%
    \futurelet\first@char\gobble@till@marker#1\endmarker
    \ifcat\noexpand\first@char A\else
      \def\first@char{}%
    \fi
    \macc@nested@a\relax111{\first@char}%
  \fi
  \endgroup
}
\makeatother

\newcommand\rsmraise[1]{%
  \ifx#1\displaystyle .8\else
    \ifx#1\textstyle .8\else
      \ifx#1\scriptstyle .6\else
        .45%
      \fi
    \fi
  \fi}

\title{Quantum Howe duality and invariant polynomials}

\author{Vyacheslav Futorny, Libor Křižka, Jian Zhang}

\AtEndDocument{\bigskip{\footnotesize%
  (V.\,Futorny) \textsc{Instituto de Matemática e Estatística, Universidade de Sa\~{o} Paulo, Caixa Postal 66281,\\ Sa\~{o} Paulo, CEP 05315-970, Brasil} \par
  \textit{E-mail address}: \texttt{futorny@ime.usp.br} \par
  \addvspace{\medskipamount}
  (L.\,Křižka) \textsc{Instituto de Matemática e Estatística, Universidade de Sa\~{o} Paulo, Caixa Postal 66281,\\ Sa\~{o} Paulo, CEP 05315-970, Brasil} \par
  \textit{E-mail address}: \texttt{krizka.libor@gmail.com} \par
  \addvspace{\medskipamount}
  (J.\,Zhang) \textsc{Instituto de Matemática e Estatística, Universidade de Sa\~{o} Paulo, Caixa Postal 66281,\\ Sa\~{o} Paulo, CEP 05315-970, Brasil} \par
  \textit{E-mail address}: \texttt{zhang@ime.usp.br} \par
}}

\date{}

\begin{document}

\maketitle

\begin{abstract}

We construct two examples of $q$-deformed classical Howe dual pairs  $(\mathfrak{sl}(2,\mathbb{C}),\mathfrak{so}(3,\mathbb{C}))$ and $(\mathfrak{sl}(2,\mathbb{C}),\mathfrak{sl}(n,\mathbb{C}))$. Moreover, we obtain a noncommutative version of the first fundamental theorem of classical invariant theory. Our approach to these duality differs from \cite{Lehrer-Zhang-Zhang2011} and \cite{Noumi-Umeda-Wakayama1996}.
Furthermore, we solve the tensor product decomposition problem for Verma modules over $U_q(\mathfrak{sl}(2,\mathbb{C}))$ provided $q$ is not a root of unity.



\end{abstract}

\thispagestyle{empty}

\tableofcontents


\section*{Introduction}
\addcontentsline{toc}{section}{Introduction}

Classical Howe dualities provide a representation theoretical framework for classical invariant theory, in particular for the first and second fundamental theorem \cite{Howe1989a}.
On the other hand, quantum Howe dualities attracted attention from the very beginning of quantum group theory. The first studied case of a $q$-deformed Howe dual pair $(U_{q^2}(\mfrak{sl}(2,\C)),U_q(\mfrak{so}(n,\C)))$ was described in \cite{Noumi-Umeda-Wakayama1996}. A quantum version of the Howe duality for a pair
of quantized universal enveloping algebras of general linear algebras $U_q(\mfrak{gl}(n,\C))$ and  $U_q(\mfrak{gl}(m,\C))$ at generic $q$ is given in \cite{Zhang2003}. Furthermore, a noncommutative version of the first fundamental theorem for quantum groups of classical Lie algebras was established in \cite{Lehrer-Zhang-Zhang2011}. For each quantum group associated with a classical Lie algebra a module is constructed which has a structure of a noncommutative associative algebra preserved by the quantum group. The subspace of invariants is finitely generated subalgebra, see \cite[Theorem 6.10]{Lehrer-Zhang-Zhang2011}.

In the current paper we present other examples of $q$-deformed Howe dual pairs. Let us briefly summarize the content of our article. In Section \ref{section-weyl}, we recall the definition of the quantum Weyl algebra of a complex vector space and the quantum Fourier transform. In Section \ref{section-strange}, we focus on the quantum group $U_q(\mfrak{sl}(2,\C))$ and its action on the coordinate algebra $\C_q[V]$ of the $3$-dimensional   representation $V$ of $U_q(\mfrak{sl}(2,\C))$ and introduce a $q$-analogue of invariant polynomials and harmonic polynomials.
We establish a new formalism for the Howe duality $(U_{q^2}(\mfrak{sl}(2,\C)),U_q(\mfrak{sl}(2,\C)))$ (see Theorem \ref{thm:Fischer decomposition U_q(sl(2,C))}),  which is a special case of the duality obtained in \cite{Noumi-Umeda-Wakayama1996}.  The choice of the $3$-dimensional   representation $V$ of $U_q(\mfrak{sl}(2,\C))$  is based on the fact that $V$ is only  flat $q$-deformation with nontrivial invariants.
  We describe the center ${\rm Z}(\C_q[V])$ of the coordinate algebra $\C_q[V]$ (see Theorem \ref{thm:invariant polynomials and center}).

  In Section \ref{section-third}, we discuss the Howe duality $(U_q(\mfrak{sl}(2,\C)),U_q(\mfrak{sl}(n,\C)))$. The main results are Theorem \ref{thm:Fischer decomposition U_q(sl(n,C))} and Theorem \ref{thm:invariant polynomials and center U_q(sl(n,C))}. Let us note that our treatment of this duality differs from \cite{Lehrer-Zhang-Zhang2011}. In Section \ref{section-verma}, we use the geometric realization of Verma modules over $U_q(\mfrak{sl}(2,\C))$ to solve the tensor product decomposition problem for Verma modules over $U_q(\mfrak{sl}(2,\C))$ (Theorem \ref{thm:tensor product Verma quantum}). Special cases of Theorem \ref{thm:tensor product Verma quantum} were obtained previously in \cite[Section 3.5]{Chari-Jakelic-Moura2005} for  $q^{2\lambda}, q^{2\mu} \in q^{2\N_0}$ (see also \cite{Creath-Jakelic2018}) and in \cite[Section 2]{Frenkel-Zeitlin2010}, \cite[Section 5]{Huang2017} for $q^{2(\lambda+\mu)}, q^{2\lambda}, q^{2\mu} \notin q^{2\N_0}$.

Throughout the article, unless otherwise stated, we assume that $q \in \C^\times$ is \emph{not a root of unity}. Further, we use the standard notation $\N$ and $\N_0$ for the set of natural numbers and the set of natural numbers together with zero, respectively.


\section{Quantum Weyl algebra and quantum Fourier transform}
\label{section-weyl}

For $q \in \C^\times$ satisfying $q \neq \pm 1$ and $a \in \C$, the $q$-number $[a]_q$ is defined by
\begin{align}
  [a]_q = {q^a-q^{-a} \over q-q^{-1}}.
\end{align}
If $n \in \N_0$, then we introduce the $q$-factorial $[n]_q!$ by
\begin{align}
  [n]_q!= \prod_{k=1}^n [k]_q.
\end{align}
The $q$-binomial coefficients are defined by the formula
\begin{align}
  {n \brack k}_q = {[n]_q! \over [k]_q![n-k]_q!},
\end{align}
where $n, k \in \N_0$ and $n \geq k$.
\medskip

Let us consider an associative $\C$-algebra $A$. Let $\sigma \colon A \rarr A$ be a $\C$-algebra automorphism. Then a twisted derivation of $A$ relative to $\sigma$ is a linear mapping $D \colon A \rarr A$ satisfying
\begin{align}
  D(ab) = D(a)\sigma(b) + \sigma^{-1}(a) D(b)
\end{align}
for all $a,b \in A$. An element $a \in A$ induces an inner twisted derivation $\ad_\sigma\! a$ relative to $\sigma$ defined by the formula
\begin{align}
  (\ad_\sigma\! a)(b)= a \sigma(b) - \sigma^{-1}(b)a
\end{align}
for all $a,b \in A$. Let us note that also $D_\sigma = \sigma - \sigma^{-1}$ is a twisted derivation of $A$ relative to $\sigma$.
\medskip

\lemma{Let $D$ be a twisted derivation of $A$ relative to $\sigma$. Then we have
\begin{align}
\begin{aligned}
  \sigma \circ \lambda_a &= \lambda_{\sigma(a)} \circ \sigma, \qquad &    D \circ \lambda_a - \lambda_{\sigma^{-1}(a)} \circ D &= \lambda_{D(a)} \circ \sigma, \\
  \sigma \circ \rho_a &= \rho_{\sigma(a)} \circ \sigma, &
   D \circ \rho_a - \rho_{\sigma(a)} \circ D &= \rho_{D(a)} \circ \sigma^{-1}
\end{aligned}
\end{align}
for all $a \in A$, where $\lambda_a$ and $\rho_a$ denote the left and the right multiplications by $a \in A$, respectively.}

\proof{We have
\begin{gather*}
  (\sigma \circ \lambda_a)(b)= \sigma(ab)= \sigma(a)\sigma(b) = (\lambda_{\sigma(a)} \circ \sigma)(b) \\
  (D \circ \lambda_a)(b)= D(ab) = D(a)\sigma(b)+\sigma^{-1}(a)D(b) = (\lambda_{D(a)}\circ \sigma + \lambda_{\sigma^{-1}(a)} \circ D)(b)
\end{gather*}
for all $a,b \in A$.}

Let $V$ be a finite-dimensional complex vector space and let $\C[V]$ be the $\C$-algebra of polynomial functions on $V$. Further, let $\{x_1,x_2,\dots,x_n\}$ be the linear coordinate functions on $V$ with respect to a basis $\{e_1,e_2,\dots,e_n\}$ of $V$. Then there exists a canonical isomorphism of $\C$-algebras $\C[V]$ and $\C[x_1,x_2,\dots,x_n]$.

Let $q \in \C^\times$ satisfies $q \neq \pm 1$. We define a $\C$-algebra automorphism $\gamma_{q,x_i}$ of $\C[V]$ by
\begin{align}
  \gamma_{q,x_i} = q^{x_i \partial_{x_i}}
\end{align}
and a twisted derivation $\partial_{q,x_i}$ of $\C[V]$ relative to $\gamma_{q,x_i}$ through
\begin{align}
  \partial_{q,x_i} = {1 \over x_i} {q^{x_i\partial_{x_i}} - q^{-x_i\partial_{x_i}} \over q - q^{-1}}
\end{align}
for $i=1,2,\dots,n$.
\medskip

\lemma{Let $q \in \C$ satisfies $q \neq \pm 1$. Further, let $D$ be a twisted derivation of $\C[V]$ relative to $\gamma_{q,x_i}$ for some $i=1,2,\dots,n$. Then we have
\begin{align}
  D =   f_i \partial_{q,x_i},
\end{align}
where $f_i \in \C[V]$.}

\proof{For $j=1,2,\dots,n$ satisfying $j \neq i$, we have
\begin{align*}
  D(x_ix_j) &= D(x_i)\gamma_{q,x_i}(x_j)+\gamma_{q,x_i}^{-1}(x_i)D(x_j) =x_jD(x_i)+q^{-1}x_iD(x_j), \\
  D(x_jx_i) &= D(x_j)\gamma_{q,x_i}(x_i)+\gamma_{q,x_i}^{-1}(x_j)D(x_i) =qx_iD(x_j)+x_jD(x_i),
\end{align*}
which implies that $D(x_j)=0$ for all $j=1,2,\dots,n$ such that $j \neq i$. If we set $f_i=D(x_i)$, then we get
\begin{align*}
  (D-f_i\partial_{q,x_i})(x_j)=0
\end{align*}
for all $j=1,2,\dots,n$, which gives us $D=f_i \partial_{q,x_i}$.}

Let $q \in \C^\times$ satisfies $q \neq \pm 1$. Then based on the previous lemma, we define the quantum Weyl algebra $\eus{A}^q_V$ of the complex vector space $V$ as an associative $\C$-subalgebra of $\End \C[V]$ generated by $x_i$, $\partial_{q,x_i}$ and $\gamma_{q,x_i}^{\pm 1}$ for $i=1,2,\dots,n$. Let us note that the definition of $\eus{A}^q_V$ depends on the choice of a basis $\{e_1,e_2,\dots,e_n\}$ of $V$.

Moreover, we have the following nontrivial relations
\begin{align}
 \gamma_{q,x_i} x_i  = q x_i \gamma_{q,x_i}, \qquad \gamma_{q,x_i} \partial_{q,x_i} = q^{-1} \partial_{q,x_i}\gamma_{q,x_i}
\end{align}
and
\begin{align}
  \partial_{q,x_i} x_i - q x_i \partial_{q,x_i}=\gamma_{q,x_i}^{-1}, \qquad   \partial_{q,x_i} x_i - q^{-1} x_i \partial_{q,x_i}=\gamma_{q,x_i}
\end{align}
for $i=1,2,\dots,n$.

Furthermore, let $\{y_1,y_2,\dots,y_n\}$ be the dual linear coordinate functions on $V^*$. Then there is a canonical isomorphism
\begin{align}
  \mcal{F} \colon \eus{A}^q_V \rarr \eus{A}^q_{V^*}
\end{align}
of associative $\C$-algebras given by
\begin{align}
  \mcal{F}(x_i) = -\partial_{q,y_i}, \qquad \mcal{F}(\partial_{q,x_i})=y_i, \qquad  \mcal{F}(\gamma_{q,x_i}) = q^{-1} \gamma_{q,y_i}^{-1}
\end{align}
for $i=1,2,\dots,n$. We call $\mcal{F} \colon \eus{A}^q_V \rarr \eus{A}^q_{V^*}$ the quantum Fourier transform.


\section{Representations of $U_q(\mfrak{sl}(2,\C))$ and its Howe duality}
\label{section-strange}

In this section we will focus more concretely on the simplest possible quantum group $U_q(\mfrak{sl}(2,\C))$.


\subsection{The quantum group $U_q(\mfrak{sl}(2,\C))$}

Let us consider the complex simple Lie algebra $\mfrak{sl}(2,\C)$. We denote by
\begin{align}
  e =\begin{pmatrix}
    0 & 1 \\
    0 & 0
  \end{pmatrix}\!, \qquad
  h=\begin{pmatrix}
    1 & 0 \\
    0 & -1
  \end{pmatrix}\!, \qquad
  f=\begin{pmatrix}
    0 & 0 \\
    1 & 0
  \end{pmatrix}
\end{align}
the standard basis of $\mfrak{sl}(2,\C)$. The vector subspace $\mfrak{h}=\C h$ is a Cartan subalgebra of $\mfrak{sl}(2,\C)$. Let us define $\alpha \in \mfrak{h}^*$ by $\alpha(h)=2$. Then the root system of $\mfrak{sl}(2,\C)$ with respect to $\mfrak{h}$ is $\Delta=\{\pm\alpha\}$ and a positive root system in $\Delta$ is $\Delta^+=\{\alpha\}$. The fundamental weight is given by $\omega={1 \over 2}\alpha$. The standard Borel subalgebra $\mfrak{b}$ of $\mfrak{sl}(2,\C)$ is defined by $\mfrak{b}=\C h \oplus \C e$ with the nilradical $\mfrak{n}=\C e$ and the opposite nilradical $\widebar{\mfrak{n}} = \C f$.

Let $q \in \C^\times$ satisfies $q \neq \pm 1$. The quantum group $U_q(\mfrak{sl}(2,\C))$, see \cite{Jimbo1986}, is a unital associative $\C$-algebra generated by $E,F,K,K^{-1}$ subject to the relations
\begin{gather}
\begin{gathered}
KK^{-1}=1, \qquad K^{-1}K=1, \\
  KEK^{-1}=q^2E, \qquad [E,F]= {K-K^{-1} \over q-q^{-1}}, \qquad KFK^{-1}=q^{-2}F.
\end{gathered}
\end{gather}
There is a unique Hopf algebra structure on $U_q(\mfrak{sl}(2,\C))$ with the coproduct $\Delta \colon U_q(\mfrak{sl}(2,\C)) \rarr U_q(\mfrak{sl}(2,\C)) \otimes_\C U_q(\mfrak{sl}(2,\C))$, the counit $\veps \colon U_q(\mfrak{sl}(2,\C)) \rarr \C$ and the antipode $S \colon U_q(\mfrak{sl}(2,\C)) \rarr U_q(\mfrak{sl}(2,\C))$ given by
\begin{gather} \label{eq:Hopf algebra structure U_q(sl(2,C))}
\begin{gathered}
 \Delta(E)=E \otimes K + 1 \otimes E,\quad \Delta(K)=K \otimes K,\quad \Delta(F)=F \otimes 1 + K^{-1} \otimes F, \\
 \veps(E)=0, \qquad \qquad \veps(K) =1, \qquad \qquad \veps(F)=0, \\
 S(E)=-EK^{-1}, \qquad S(K)=K^{-1}, \qquad S(F)=-KF.
\end{gathered}
\end{gather}
\medskip

\lemma{\label{lem:commutator U_q(sl(2,C))}
Let $q \in \C^\times$ satisfies $q \neq \pm 1$. Then in the quantum group $U_q(\mfrak{sl}(2,\C))$ we have
\begin{align} \label{eq:commutator U_q(sl(2,C))}
\begin{aligned}[]
  [E^s,F]&= [s]_q E^{s-1} {q^{s-1}K - q^{-s+1}K^{-1} \over q-q^{-1}},\\
  [E,F^s]&= [s]_q F^{s-1} {q^{-s+1}K - q^{s-1}K^{-1} \over q-q^{-1}}
\end{aligned}
\end{align}
for all $s \in \N_0$.}

\vspace{-2mm}



\subsection{Quantum invariant polynomials}

Let $V$ be a finite-dimensional representation of a complex semisimple Lie algebra $\mfrak{g}$. Then we have the induced representation of $\mfrak{g}$ on the $\C$-algebra $\C[V]$ of polynomial functions on $V$. In \cite{Lusztig1988} and \cite{Rosso1988} it was proved that $V$ can be $q$-deformed into a finite-dimensional representation of the quantum group $U_q(\mfrak{g})$. Hence, it arises a natural question of a $q$-deformation of $\C[V]$. In the next, we shall describe this $q$-deformation of $\C[V]$ in one particular example, which offers a good insight into a general construction. Let us note that this construction works for any finite-dimensional representation $V$ of a semisimple Lie algebra $\mfrak{g}$ not only of $\mfrak{sl}(2,\C)$.  The $\C$-algebra $\C_q[V]$ is usually called the coordinate algebra of the quantum vector space $V$ introduced in \cite{Reshetikhin-Takhtajan-Faddeev1990}.

\medskip

Let $V^*=\C^3{}^*$ be a $3$-dimensional simple $U_q(\mfrak{sl}(2,\C))$-module defined by
\begin{align} \label{eq:adjoint representation U_q(sl(2,C))}
  \rho_q(E)=\begin{pmatrix}
    0 & 0 & 0 \\
    [2]_q & 0 & 0\\
    0 & -1 & 0
  \end{pmatrix}\!, \qquad
  \rho_q(K)=\begin{pmatrix}
    q^{-2} & 0 & 0 \\
    0 & 1 & 0 \\
    0 & 0 & q^2
  \end{pmatrix}\!,
  \qquad
  \rho_q(F)=\begin{pmatrix}
    0 & 1 & 0 \\
    0 & 0 & -[2]_q \\
    0 & 0 &  0
  \end{pmatrix}
\end{align}
with respect to the canonical basis $(x,z,y)$ of $\C^3{}^*$. Since $V^* \otimes_\C V^* \simeq L_q(4\omega) \oplus L_q(2\omega) \oplus L_q(0)$, where $L_q(\lambda\omega)$ is the simple $U_q(\mfrak{sl}(2,\C))$-module with highest weight $q^{\lambda \omega}$ for $\lambda \in \C$, we define
\begin{align}
 \C_q[V] \simeq S_q(V^*) = T(V^*)/I_q,
\end{align}
where $I_q$ is the two-sided ideal of the tensor algebra $T(V^*)$ generated by $L_q(2\omega)=\langle x \otimes z-q^2z \otimes x,x \otimes y-y \otimes x+q(q^2-q^{-2})z \otimes z,y\otimes z-q^{-2}z \otimes y \rangle $, which gives us
\begin{align}
  \C_q[V] \simeq \C_q[x,y,z] = \C\langle x,y,z \rangle /(xz-q^2zx,xy-yx+q(q^2-q^{-2})z^2,yz-q^{-2}zy). \label{eq:adjoint representation quantum ring}
\end{align}
Moreover, since the two-sided ideal $I_q$ is a $U_q(\mfrak{sl}(2,\C))$-submodule of $T(V^*)$, we obtain that $\C_q[ V]$ is a $U_q(\mfrak{sl}(2,\C))$-module. Taking the limit $q \rarr 1$, we get $I_q \rarr I$, and hence $\C_q[x,y,z] \rarr \C[x,y,z]$.

As the $\C$-algebra $\C_q[x,y,z]$ has a basis $\{x^az^cy^b;\, a,b,c \in \N_0\}$, we can find a family of isomorphisms $\varphi_q \colon \C[x,y,z] \rarr \C_q[x,y,z]$ of vector spaces such that $\varphi_q \rarr \id$ for $q \rarr 1$. Let us define $\varphi_q \colon \C[x,y,z] \rarr \C_q[x,y,z]$ by
\begin{align}
  \varphi_q(x^ay^bz^c)=x^az^cy^b \label{eq:isomomorphism polynomials}
\end{align}
for all $a,b,c \in \N_0$. Then the corresponding $U_q(\mfrak{sl}(2,\C))$-module structure on $\C[x,y,z]$ is given through the homomorphism
\begin{align}
\sigma_q \colon U_q(\mfrak{sl}(2,\C)) \rarr \eus{A}^q_V
\end{align}
of associative $\C$-algebras, where $\eus{A}^q_V$ is the quantum Weyl algebra of the vector space $V$, defined by
\begin{align}
  \sigma_q(a) = \varphi_q^{-1}\! \circ \rho_q(a) \circ \varphi_q
\end{align}
for all $a \in U_q(\mfrak{sl}(2,\C))$.  Let us note that a different choice of the isomorphism $\varphi_q$ leads only to a distinct realization of the $U_q(\mfrak{sl}(2,\C))$-module $\C_q[V]$.
\medskip

\proposition{The homomorphism $\sigma_q \colon U_q(\mfrak{sl}(2,\C)) \rarr \eus{A}^q_V$ of associative $\C$-algebras is given by
\begin{align}
\begin{aligned}
  \sigma_q(E)&= [2]_qz\gamma_{q,x}^{-2} \gamma_{q,y}^2\partial_{q^2,x} -y \gamma_{q,y}^2 \gamma_{q,z}^{-1}\partial_{q,z}, \\
  \sigma_q(K)&=\gamma_{q,x}^{-2} \gamma_{q,y}^2, \\
  \sigma_q(F)&=-[2]_q z\gamma_{q,x}^2 \gamma_{q,y}^{-2} \partial_{q^2,y} + x \gamma_{q,x}^2 \gamma_{q,z}^{-1}\partial_{q,z}.
\end{aligned}
\end{align}
}

\proof{The proof is a straightforward computation. Using \eqref{eq:adjoint representation U_q(sl(2,C))} and \eqref{eq:Hopf algebra structure U_q(sl(2,C))}, we may write
\begin{align*}
  \rho_q(K)(x^az^cy^b)=q^{-2a+2b}x^az^cy^b,
\end{align*}
\begin{align*}
  \rho_q(E)(x^az^cy^b)&=\sum_{k=1}^a q^{2b-2a+2k}[2]_q x^{k-1}zx^{a-k}z^cy^b -\sum_{k=1}^c q^{2b}x^a z^{k-1}yz^{c-k}y^b  \\
  &=\sum_{k=1}^a q^{2b-4a+4k}[2]_q x^{a-1}z^{c+1}y^b -\sum_{k=1}^c q^{2b-2c+2k}x^a z^{c-1}y^{b+1} \\&=[2]_qq^{2b}{1-q^{-4a} \over 1-q^{-4}}\, x^{a-1}z^{c+1}y^b - q^{2b} {1-q^{-2c} \over 1-q^{-2}}\,x^a z^{c-1}y^{b+1} \\
  &= [2]_q q^{2b-2a+2} [a]_{q^2} x^{a-1}z^{c+1}y^b - q^{2b-c+1}[c]_q x^a z^{c-1}y^{b+1}
\end{align*}
and
\begin{align*}
  \rho_q(F)(x^az^cy^b)&=- \sum_{k=1}^b q^{2a-2k+2}[2]_q x^az^cy^{k-1}zy^{b-k}+\sum_{k=1}^c q^{2a}x^a z^{k-1}xz^{c-k}y^b    \\
  &=- \sum_{k=1}^b q^{2a-4k+4}[2]_q x^az^{c+1}y^{b-1} +\sum_{k=1}^c q^{2a-2k+2}x^{a+1} z^{c-1}y^b \\&=-[2]_qq^{2a}{1-q^{-4b} \over 1-q^{-4}}\, x^az^{c+1}y^{b-1} + q^{2a} {1-q^{-2c} \over 1-q^{-2}}\,x^{a+1} z^{c-1}y^b \\
  &= -[2]_q q^{2a-2b+2} [b]_{q^2} x^az^{c+1}y^{b-1} + q^{2a-c+1}[c]_q x^{a+1} z^{c-1}y^b
\end{align*}
for all $a,b,c \in \N_0$. Therefore, we obtain
\begin{align*}
  \sigma_q(E)(x^ay^bz^c)&=[2]_q q^{2b-2a+2} [a]_{q^2} x^{a-1}y^bz^{c+1} - q^{2b-c+1}[c]_q x^a y^{b+1}z^{c-1} \\
  &=[2]_qz\gamma_{q,x}^{-2}\gamma_{q,y}^2 \partial_{q^2,x}(x^ay^bz^c)-y\gamma_{q,y}^2 \gamma_{q,z}^{-1} \partial_{q,z}(x^ay^bz^c), \\[1mm]
  \sigma_q(K)(x^ay^bz^c)&=q^{-2a+2b}x^ay^bz^c=\gamma_{q,x}^{-2}\gamma_{q,y}^2(x^ay^bz^c), \\[1mm]
   \sigma_q(F)(x^ay^bz^c)&=-[2]_q q^{2a-2b+2} [b]_{q^2} x^ay^{b-1}z^{c+1}+ q^{2a-c+1}[c]_q x^{a+1}y^bz^{c-1} \\
   &=-[2]_qz\gamma_{q,x}^2 \gamma_{q,y}^{-2} \partial_{q^2,y}(x^ay^bz^c)+x\gamma_{q,x}^2 \gamma_{q,z}^{-1} \partial_{q,z}(x^ay^bz^c)
\end{align*}
for all $a,b,c \in \N_0$.}
\vspace{-2mm}


\subsection{The Howe duality $(U_{q^2}(\mfrak{sl}(2,\C)),U_q(\mfrak{sl}(2,\C)))$}

In this section we present a quantum analogue of the classical Howe duality $(\mfrak{sl}(2,\C), \mfrak{so}(3,\C))$. We describe this duality by means of quantum differential operators on a complex vector space, which provides a general tool for investigating other quantum Howe dualities.
\medskip

The polynomials belonging to
\begin{align}
  \C_q[V]^{U_q(\mfrak{sl}(2,\C))} = \{f \in \C_q[V];\, \rho_q(a)f = \veps(a)f\ \text{for all $a \in U_q(\mfrak{sl}(2,\C))$}\}
\end{align}
are called invariant polynomials. Moreover, the set $\C_q[V]^{U_q(\mfrak{sl}(2,\C))}$ of all invariant polynomials is a $\C$-subalgebra of $\C_q[V]$. In the next, we determine this algebra of invariant polynomials using a $q$-analogue of the Fischer decomposition.
\medskip

Let us introduce the differential operators
\begin{align}
\begin{gathered}
  P_q= q^{-1}xy\gamma_{q,z}^{-2}+q^2z^2, \qquad Q_q=\partial_{q^2,x} \partial_{q^2,y}+{\textstyle {1 \over [2]_q^2}}\, \partial_{q,z}^2 \gamma_{q,x}^2 \gamma_{q,y}^2 \\
  E_q= \gamma_{q,x} \gamma_{q,y} \gamma_{q,z}
\end{gathered}
\end{align}
in the quantum Weyl algebra $\eus{A}^q_V$.
\medskip

\proposition{\label{prop:invariant operators and dual group}
The differential operators $P_q$, $Q_q$ and $E_q$ are $U_q(\mfrak{sl}(2,\C))$-invariant, i.e.\ they commute with $\sigma_q(a)$ for all $a \in U_q(\mfrak{sl}(2,\C))$. Moreover, the mapping $\pi_q \colon U_{q^2}(\mfrak{sl}(2,\C)) \rarr \eus{A}^q_V$ uniquely determined by
\begin{align}
  \pi_q(E) = P_q, \qquad \pi_q(K)=q^3 E_q^2, \qquad \pi_q(F)=-Q_q
\end{align}
gives rise to a homomorphism of associative $\C$-algebras.}

\proof{We may write
\begin{align*}
  [\sigma_q(E),P_q]&=[[2]_qz\gamma_{q,x}^{-2} \gamma_{q,y}^2\partial_{q^2,x} -y \gamma_{q,y}^2 \gamma_{q,z}^{-1}\partial_{q,z},q^{-1}xy\gamma_{q,z}^{-2}+q^2z^2]\\
  &=[[2]_qz\gamma_{q,x}^{-2} \gamma_{q,y}^2\partial_{q^2,x},q^{-1}xy\gamma_{q,z}^{-2}]-[y \gamma_{q,y}^2 \gamma_{q,z}^{-1}\partial_{q,z},q^2z^2]  \\
  &=[2]_qqy\gamma_{q,x}^{-2} \gamma_{q,y}^2[z\partial_{q^2,x},x\gamma_{q,z}^{-2}]-q^2y \gamma_{q,y}^2[\gamma_{q,z}^{-1}\partial_{q,z},z^2]  \\
  &=[2]_qqy\gamma_{q,x}^{-2} \gamma_{q,y}^2z\gamma_{q,x}^2\gamma_{q,z}^{-2}-q^2y \gamma_{q,y}^2(1+q^{-2})z\gamma_{q,z}^{-2}=0, \\
  [\sigma_q(K),P_q]&=[\gamma_{q,x}^{-2} \gamma_{q,y}^2,q^{-1}xy\gamma_{q,z}^{-2}+q^2z^2]=0, \\
  [\sigma_q(F),P_q]&=[-[2]_q z\gamma_{q,x}^2 \gamma_{q,y}^{-2} \partial_{q^2,y} + x \gamma_{q,x}^2 \gamma_{q,z}^{-1}\partial_{q,z},q^{-1}xy\gamma_{q,z}^{-2}+q^2z^2] \\
  &=-[[2]_q z\gamma_{q,x}^2 \gamma_{q,y}^{-2} \partial_{q^2,y},q^{-1}xy\gamma_{q,z}^{-2}]+ [x \gamma_{q,x}^2 \gamma_{q,z}^{-1}\partial_{q,z},q^2z^2] \\
  &=-[2]_qqx\gamma_{q,x}^2 \gamma_{q,y}^{-2}[z\partial_{q^2,y},y\gamma_{q,z}^{-2}]+ q^2x \gamma_{q,x}^2[\gamma_{q,z}^{-1}\partial_{q,z},z^2] \\
  &= -[2]_qqx\gamma_{q,x}^2 \gamma_{q,y}^{-2}z\gamma_{q,y}^2\gamma_{q,z}^{-2} + q^2x \gamma_{q,x}^2(1+q^{-2})z\gamma_{q,z}^{-2}=0
\end{align*}
and
\begin{align*}
  [\sigma_q(E),Q_q]&=[[2]_qz\gamma_{q,x}^{-2} \gamma_{q,y}^2\partial_{q^2,x} -y \gamma_{q,y}^2 \gamma_{q,z}^{-1}\partial_{q,z},\partial_{q^2,x} \partial_{q^2,y}+{\textstyle {1 \over [2]_q^2}}\, \partial_{q,z}^2 \gamma_{q,x}^2 \gamma_{q,y}^2]\\
  &=[[2]_qz\gamma_{q,x}^{-2} \gamma_{q,y}^2\partial_{q^2,x}, {\textstyle {1 \over [2]_q^2}}\, \partial_{q,z}^2 \gamma_{q,x}^2 \gamma_{q,y}^2] - [y \gamma_{q,y}^2 \gamma_{q,z}^{-1}\partial_{q,z},\partial_{q^2,x} \partial_{q^2,y}] \\
  &={\textstyle {1 \over [2]_q}}\, \gamma_{q,y}^4[z\gamma_{q,x}^{-2} \partial_{q^2,x}, \partial_{q,z}^2 \gamma_{q,x}^2] - \gamma_{q,z}^{-1}\partial_{q,z}\partial_{q^2,x}[y \gamma_{q,y}^2 , \partial_{q^2,y}] \\
  &={\textstyle {1 \over [2]_q}}\, \gamma_{q,y}^4(-[2]_q\gamma_{q,z}^{-1}\partial_{q,z} \partial_{q^2,x}) + \gamma_{q,z}^{-1}\partial_{q,z}\partial_{q^2,x} \gamma_{q,y}^4=0, \\
  [\sigma_q(K),Q_q]&=[\gamma_{q,x}^{-2} \gamma_{q,y}^2,\partial_{q^2,x} \partial_{q^2,y}+{\textstyle {1 \over [2]_q^2}}\, \partial_{q,z}^2 \gamma_{q,x}^2 \gamma_{q,y}^2 ]=0, \\
  [\sigma_q(F),Q_q]&=[-[2]_q z\gamma_{q,x}^2 \gamma_{q,y}^{-2} \partial_{q^2,y} + x \gamma_{q,x}^2 \gamma_{q,z}^{-1}\partial_{q,z}, \partial_{q^2,x} \partial_{q^2,y}+{\textstyle {1 \over [2]_q^2}}\, \partial_{q,z}^2 \gamma_{q,x}^2 \gamma_{q,y}^2] \\
  &=-[[2]_q z\gamma_{q,x}^2 \gamma_{q,y}^{-2} \partial_{q^2,y},{\textstyle {1 \over [2]_q^2}}\, \partial_{q,z}^2 \gamma_{q,x}^2 \gamma_{q,y}^2] +[x \gamma_{q,x}^2 \gamma_{q,z}^{-1}\partial_{q,z}, \partial_{q^2,x} \partial_{q^2,y}] \\
  &= -{\textstyle {1 \over [2]_q}}\,\gamma_{q,x}^4[z \gamma_{q,y}^{-2} \partial_{q^2,y}, \partial_{q,z}^2 \gamma_{q,y}^2] +\gamma_{q,z}^{-1}\partial_{q,z}\partial_{q^2,y}[x \gamma_{q,x}^2 , \partial_{q^2,x}] \\
  &= -{\textstyle {1 \over [2]_q}}\,\gamma_{q,x}^4(-[2]_q\gamma_{q,z}^{-1}\partial_{q,z} \partial_{q^2,y}) -\gamma_{q,z}^{-1}\partial_{q,z}\partial_{q^2,y} \gamma_{q,x}^4 =0.
\end{align*}
Furthermore, we immediately obtain
\begin{align*}
  [\sigma_q(E),E_q]=0, \qquad [\sigma_q(K),E_q]=0, \qquad [\sigma_q(F),E_q]=0.
\end{align*}
Therefore, the differential operators $P_q$, $Q_q$ and $E_q$ are $U_q(\mfrak{sl}(2,\C))$-invariant. Moreover, we have
\begin{align*}
  \pi_q(K)\pi_q(E) &=q^3 \gamma_{q,x}^2 \gamma_{q,y}^2 \gamma_{q,z}^2(q^{-1}xy\gamma_{q,z}^{-2}+q^2z^2)\\
  &=q^4(q^{-1}xy\gamma_{q,z}^{-2}+q^2z^2)q^3 \gamma_{q,x}^2 \gamma_{q,y}^2 \gamma_{q,z}^2=q^4\pi_q(E)\pi_q(K), \\[1mm]
 \pi_q(K)\pi_q(F) &=-q^3 \gamma_{q,x}^2 \gamma_{q,y}^2 \gamma_{q,z}^2\big(\partial_{q^2,x} \partial_{q^2,y}+{\textstyle {1 \over [2]_q^2}}\, \partial_{q,z}^2 \gamma_{q,x}^2 \gamma_{q,y}^2\big)\\
  &=-q^{-4}\big(\partial_{q^2,x} \partial_{q^2,y}+{\textstyle {1 \over [2]_q^2}}\, \partial_{q,z}^2 \gamma_{q,x}^2 \gamma_{q,y}^2\big)q^3 \gamma_{q,x}^2 \gamma_{q,y}^2 \gamma_{q,z}^2=q^{-4}\pi_q(F)\pi_q(K)
\end{align*}
and
\begin{align*}
  [\pi_q(E), \pi_q(F)]&=[\partial_{q^2,x} \partial_{q^2,y}+{\textstyle {1 \over [2]_q^2}}\, \partial_{q,z}^2 \gamma_{q,x}^2 \gamma_{q,y}^2, q^{-1}xy\gamma_{q,z}^{-2}+q^2z^2] \\
  &=[\partial_{q^2,x} \partial_{q^2,y},q^{-1}xy\gamma_{q,z}^{-2}]+{\textstyle {1 \over [2]_q^2}}\, [\partial_{q,z}^2 \gamma_{q,x}^2 \gamma_{q,y}^2, q^2z^2] \\
  &=q^{-1}\gamma_{q,z}^{-2}[\partial_{q^2,x}\partial_{q^2,y},xy]+ {\textstyle {1 \over [2]_q^2}}\,q^2 \gamma_{q,x}^2 \gamma_{q,y}^2[\partial_{q,z}^2,z^2] \\
  &=q^{-1}\gamma_{q,z}^{-2}{q^2 \gamma_{q,x}^2\gamma_{q,y}^2- q^{-2}\gamma_{q,x}^{-2}\gamma_{q,y}^{-2} \over q^2-q^{-2}} +  {1 \over [2]_q^2}\,q^2 \gamma_{q,x}^2 \gamma_{q,y}^2 {[2]_q(q\gamma_{q,z}^2-q^{-1}\gamma_{q,z}^{-2}) \over q-q^{-1}} \\
  &={q^3 \gamma_{q,x}^2 \gamma_{q,y}^2 \gamma_{q,z}^2 - q^{-3}\gamma_{q,x}^{-2} \gamma_{q,y}^{-2} \gamma_{q,z}^{-2}\over q^2-q^{-2}} = {\pi_q(K) - \pi_q(K)^{-1} \over q^2-q^{-2}},
\end{align*}
which gives rise to a homomorphism $\pi_q\colon U_{q^2}(\mfrak{sl}(2,\C)) \rarr \eus{A}^q_V$ of associative $\C$-algebras.}

\proposition{\label{prop:orthogonal decompostion U_q(sl(2,C))}
We have a decompostion
\begin{align}
  \C[x,y,z]= \bigoplus_{j=0}^\infty \bigoplus_{a=0}^\infty P_q^j H_{q,a},
\end{align}
where
\begin{align}
  H_{q,a}=\{f \in \C[x,y,z];\, Q_qf=0,\, E_qf=q^af\}
\end{align}
is the simple finite-dimensional highest weight $U_q(\mfrak{sl}(2,\C))$-module with highest weight $q^{2a\omega}$ generated by the vector $y^a \in \C[x,y,z]$ with
\begin{align}
\dim H_{q,a} = 2a+1
\end{align}
for all $a \in \N_0$.}

\proof{First of all, we prove a decomposition
\begin{align}
  \C[x,y,z]_a = H_{q,a} \oplus P_q(\C[x,y,z]_{a-2}) \label{eq:direct sum U_q(sl(2,C))}
\end{align}
for all $a \in \N_0$, where $\C[x,y,z]_a$ denotes the eigenspace of $E_q$ on $\C[x,y,z]$ with eigenvalue $q^a$. Let us suppose that $f \in  H_{q,a} \cap P_q(\C[x,y,z]_{a-2})$ and let us take the maximal integer $s \in \N_0$ such that $f = P_q^sg$ with $g \neq 0$. Using \eqref{eq:commutator U_q(sl(2,C))} and Proposition \ref{prop:invariant operators and dual group} we have
\begin{align*}
  0 = Q_q(f) = Q_q (P_q^s g) = [Q_q,P_q^s](g)+P_q^sQ_q(g)=P_q^sQ_q(g)+[s]_{q^2}[2a-2s+1]_{q^2}P_q^{s-1}(g),
\end{align*}
where we used the fact that $\pi_q(E)=P_q$, $\pi_q(K)=q^3E_q^2$ and $\pi_q(F)=-Q_q$. As $P_q$ is an injective differential operator, we obtain
\begin{align*}
  0=P_qQ_q(g)+[s]_{q_2}[2a-2s+1]_{q^2}g.
\end{align*}
Since the coefficient $[s]_{q^2}[2a-2s+1]_{q^2}$ is nonzero, we may write $f=P_q^{s+1}h$ with $h \neq 0$, which is in a contradiction with the maximality of the integer $s$. Therefore, we have $H_{q,a} \cap P_q(\C[x,y,z]_{a-2}) = \{0\}$. Further, since the mapping $Q_q \colon \C[x,y,z]_a \rarr \C[x,y,z]_{a-2}$ is injective on $P_q(\C[x,y,z]_{a-2}) \simeq \C[x,y,z]_{a-2}$, it is also surjective for all $a \in \N_0$, which implies that $\C[x,y,z]_a$ is a direct sum of $H_{q,a}$ and $P_q(\C[x,y,z]_{a-2})$ for all $a\in \N_0$.

As a consequence of the decomposition \eqref{eq:direct sum U_q(sl(2,C))} we immediately obtain a decomposition
\begin{align}
  \C[x,y,z]_a= \bigoplus_{j=0}^{\lfloor {a \over 2} \rfloor} P_q^j H_{q,a-2j}, \label{eq:full decomposition U_q(sl(2,C))}
\end{align}
for all $a \in \N_0$. By Proposition \ref{prop:invariant operators and dual group} the differential operators $P_q$, $E_q$ and $Q_q$ are $U_q(\mfrak{sl}(2,\C))$-invariant, therefore the subspaces $P_q^j H_{q,a}$ are $U_q(\mfrak{sl}(2,\C))$-submodules and $P_q^j H_{q,a} \simeq H_{q,a}$ as $U_q(\mfrak{sl}(2,\C))$-modules for all $j,a \in \N_0$. Finally, from \eqref{eq:direct sum U_q(sl(2,C))} we have
\begin{align*}
  \dim H_{q,a} = \dim \C[x,y,z]_a - \dim \C[x,y,z]_{a-2} = {\textstyle \binom{a+2}{2}}- {\textstyle \binom{a}{2}}=2a+1
\end{align*}
for all $a \in \N_0$. Since $Q_q(y^a)=0$, $\sigma_q(E)y^a=0$ and  $\sigma_q(K)y^a=q^{2a}y^a$, we obtain that $H_{q,a}$ contains a simple finite-dimensional highest weight $U_q(\mfrak{sl}(2,\C))$-submodule isomorphic to $L_q(2a\omega)$, which is however isomorphic to $H_{q,a}$, because it has the dimension as $H_{q,a}$ for all $a\in \N_0$.}

\lemma{\label{lem:tranformation of P_q}
The element $p_q=q^{-1}xy+q^2z^2=q^{-1}yx+q^{-2}z^2$ is central in $\C_q[x,y,z]$ and
\begin{align}
  \varphi_q \circ P_q \circ \varphi_q^{-1} = p_q, \qquad \varphi_q \circ Q_q \circ \varphi_q^{-1} = \Delta_q,
\end{align}
where
\begin{align}
  \Delta_q(x^az^cy^b)= [a]_{q^2}[b]_{q^2}x^{a-1}z^cy^{b-1}+ {1 \over [2]_q^2} [c]_q[c-1]_q q^{2a+2b} x^az^{c-2}y^b
\end{align}
for all $a,b,c \in \N_0$.}

\proof{We may write
\begin{align*}
  q^{-1}xy+q^2z^2=q^{-1}(yx-q(q^2-q^{-2})z^2)+q^2z^2=q^{-1}yx+q^{-2}z^2,
\end{align*}
where we used the relations in $\C_q[x,y,z]$. Further, have
\begin{gather*}
  (q^{-1}xy+q^2z^2)x=x(q^{-1}xy+(q^2-q^{-2})z^2)+ q^{-2}xz^2=x(q^{-1}xy+q^2z^2), \\[1mm]
  (q^{-1}xy+q^2z^2)z=z(q^{-1}xy+q^2z^2), \\[1mm]
  (q^{-1}yx+q^{-2}z^2)y= y(q^{-1}yx-(q^2-q^{-2})z^2)+ q^2yz^2=y(q^{-1}yx+q^{-2}z^2),  \\[1mm]
  (q^{-1}yx+q^{-2}z^2)z= q^{-1}yxz+q^{-2}z^3 =z(q^{-1}yx+q^{-2}z^2),
\end{gather*}
which implies that $q^{-1}xy+q^2z^2$ is a central element in $\C_q[x,y,z]$.
Finally, we obtain
\begin{align*}
  (\varphi_q \circ P_q \circ \varphi_q^{-1})(x^az^cy^b) &=\varphi_q(P_q(x^ay^bz^c))=\varphi_q(q^{-2c-1}x^{a+1}y^{b+1}z^c+q^2x^ay^bz^{c+2}) \\
  &=q^{-2c-1}x^{a+1}z^cy^{b+1}+q^2x^a z^{c+2} y^b =q^{-1} x^{a+1}yz^cy^b+q^2x^a z^{c+2} y^b \\
  &=x^a(q^{-1}xy+q^2z^2)z^cy^b=(q^{-1}xy+q^2z^2)x^az^cy^b
\intertext{and}
  (\varphi_q \circ Q_q \circ \varphi_q^{-1})(x^az^cy^b)&=\varphi_q(Q_q(x^ay^bz^c)) \\
  &= \varphi_q\big([a]_{q^2}[b]_{q^2}x^{a-1}y^{b-1}z^c+ {\textstyle {1 \over [2]_q^2}} [c]_q [c-1]_q q^{2a+2b} x^ay^bz^{c-2}\big) \\
  &=[a]_{q^2}[b]_{q^2}x^{a-1}z^cy^{b-1}+ {\textstyle {1 \over [2]_q^2}} [c]_q[c-1]_q q^{2a+2b} x^az^{c-2}y^b
\end{align*}
for all $a,b,c \in \N_0$.}

From Proposition \ref{prop:invariant operators and dual group} follows that $\C[x,y,z]$ is a $U_{q^2}(\mfrak{sl}(2,\C))$-module. Hence, also $\C_q[x,y,z]$ has a $U_{q^2}(\mfrak{sl}(2,\C))$-module structure given through the homomorphism
\begin{align}
  \tau_q\colon U_{q^2}(\mfrak{sl}(2,\C)) \rarr \End \C_q[x,y,z]
\end{align}
of associative $\C$-algebras, which is defined by
\begin{align}
  \tau_q(a)= \varphi_q \circ \pi_q(a) \circ \varphi_q^{-1}
\end{align}
for $a \in U_{q^2}(\mfrak{sl}(2,\C))$. Moreover, we have the original $U_q(\mfrak{sl}(2,\C))$-module structure on $\C_q[x,y,z]$, which commutes with the action of $U_{q^2}(\mfrak{sl}(2,\C))$ as follows from Proposition \ref{prop:invariant operators and dual group}. Therefore, we may decompose $\C_q[x,y,z]$ with respect to the action of $U_{q^2}(\mfrak{sl}(2,\C)) \otimes_\C U_q(\mfrak{sl}(2,\C))$. The corresponding decomposition is the first main result of this section.
\medskip

We denote by
\begin{align}
  \mcal{I}_q = \C[p_q]
\end{align}
the algebra of quantum invariant polynomials generated by $p_q$ and by
\begin{align}
  \mcal{H}_q = \{f \in \C_q[x,y,z];\, \Delta_qf=0\}
\end{align}
the vector space of quantum harmonic polynomials. Furthermore, we have the decomposition
\begin{align}
  \mcal{H}_q = \bigoplus_{a=0}^\infty \mcal{H}_{q,a},
\end{align}
where $\mcal{H}_{q,a}$ is the vector space of homogeneous quantum harmonic polynomials with eigenvalue $q^a$ due to $E_q$.
\medskip

\theorem{\label{thm:Fischer decomposition U_q(sl(2,C))}
We have a decomposition
\begin{align}
  \C_q[V] \simeq \mcal{I}_q \otimes_\C \mcal{H}_q = \bigoplus_{a=0}^\infty \mcal{I}_q \otimes_\C \mcal{H}_{q,a}.
\end{align}
Moreover, we have
\begin{align}
  \mcal{I}_q \otimes_\C \mcal{H}_{q,a} \simeq M_{q^2}((\textstyle{3 \over 2}+a)\omega) \otimes_\C L_q(2a\omega)
\end{align}
as $(U_{q^2}(\mfrak{sl}(2,\C)) \otimes_\C U_q(\mfrak{sl}(2,\C)))$-modules, where $M_{q^2}(\lambda\omega)$ for $\lambda \in \C$ is the Verma module with lowest weight $q^{2\lambda\omega}$ for $U_{q^2}(\mfrak{sl}(2,\C))$ and $L_q(a \omega)$ for $a \in \N_0$ is the simple finite-dimensional highest weight module with highest weight $q^{a\omega}$ for $U_q(\mfrak{sl}(2,\C))$. Furthermore, the vector $1 \otimes y^a \in \mcal{I}_q \otimes_\C \mcal{H}_{q,a}$ for $a \in \N_0$ is the (lowest, highest) weight vector with (lowest, highest) weight $(q^{(3+2a)\omega}, q^{2a\omega})$.}

\proof{From Proposition \ref{prop:orthogonal decompostion U_q(sl(2,C))} we have the decomposition
\begin{align*}
  \C[x,y,z] = \bigoplus_{j=0}^\infty \bigoplus_{a=0}^\infty P_q^j H_{q,a}.
\end{align*}
Using the isomorphism $\varphi_q \colon \C[x,y,z] \rarr \C_q[x,y,z]$ and Lemma \ref{lem:tranformation of P_q}, we may write
\begin{align*}
  \C_q[x,y,z]= \bigoplus_{j=0}^\infty \bigoplus_{a=0}^\infty\, p_q^j \varphi_q(H_{q,a}) = \bigoplus_{j=0}^\infty \bigoplus_{a=0}^\infty\, p_q^j\mcal{H}_{q,a} \simeq \bigoplus_{a=0}^\infty \mcal{I}_q \otimes_\C \mcal{H}_{q,a},
\end{align*}
since $\mcal{H}_{q,a} = \varphi_q(H_{q,a})$ for all $a\in \N_0$. Further, we have
\begin{align*}
  \sigma_q(E)(P_q^jh)=P_q^j\sigma_q(E)h,\qquad \sigma_q(K)(P_q^jh)=P_q^j\sigma_q(K)h, \qquad \sigma_q(F)(P_q^jh)&=P_q^j\sigma_q(F)h
\end{align*}
and
\begin{align*}
  \pi_q(E)(P_q^jh)&=P_q^{j+1}h, \\
  \pi_q(K)(P_q^jh)&= q^{3+2a+4j} P_q^jh,\\
  \pi_q(F)(P_q^jh)&=-Q_q(P_q^jh)=[P_q^j,Q_q](h)-P_q^jQ_q(h)=-[j]_{q^2}[\textstyle{{3 \over 2}}+a+j-1]_{q^2} P_q^{j-1}h
\end{align*}
for all $j \in \N_0$ and $h \in H_{q,a}$, where we used Lemma \ref{lem:commutator U_q(sl(2,C))} and Proposition \ref{prop:invariant operators and dual group}, which gives us that the vector subspace $\bigoplus_{j=0}^\infty P_q^j H_{q,a}$ is isomorphic to $M_{q^2}(({3 \over 2}+a)\omega) \otimes_\C L_q(2a\omega)$ as a $(U_{q^2}(\mfrak{sl}(2,\C)) \otimes_\C U_q(\mfrak{sl}(2,\C)))$-module, since Proposition \ref{prop:orthogonal decompostion U_q(sl(2,C))} implies that $H_{q,a}$ is isomorphic to $L_q(2a\omega)$ as $U_q(\mfrak{sl}(2,\C))$-module for $a \in \N_0$.}

\theorem{\label{thm:invariant polynomials and center}
The algebra of invariant quantum polynomials $\C_q[V]^{U_q(\mfrak{sl}(2,\C))}$ and the center ${\rm Z}(\C_q[V])$ are given by
\begin{align}
  {\rm Z}(\C_q[V]) = \C_q[V]^{U_q(\mfrak{sl}(2,\C))} \simeq \mcal{I}_q = \C[p_q].
\end{align}
}

\proof{By Theorem \ref{thm:Fischer decomposition U_q(sl(2,C))} we have
\begin{align*}
  \C_q[V]^{U_q(\mfrak{sl}(2,\C))} \simeq \mcal{I}_q \otimes_\C \mcal{H}_{q,0} \simeq \mcal{I}_q,
\end{align*}
since $\mcal{H}_{q,0}$ is the trivial representation of $U_q(\mfrak{sl}(2,\C))$.

Let us denote by $\lambda_a$ and $\rho_a$ the left and right multiplication by an element $a \in \C_q[x,y,z]$, respectively. Then the condition $f \in {\rm Z}(\C_q[x,y,z])$ means $(\lambda_a-\rho_a)f=0$ for all $a \in \C_q[x,y,z]$. Using the isomorphism $\varphi_q \colon \C[x,y,z] \rarr \C_q[x,y,z]$ defined by \eqref{eq:isomomorphism polynomials}, we can transfer the previous condition into $\C[x,y,z]$. Hence, we get
\begin{align*}
  (\varphi_q^{-1}\! \circ (\lambda_a-\rho_a) \circ \varphi_q)(\varphi_q^{-1}(f)) =0
\end{align*}
for all $a \in \C_q[x,y,z]$ whenever $f \in {\rm Z}(\C_q[x,y,z])$. In fact, we need to know $\varphi_q^{-1}\! \circ (\lambda_a-\rho_a) \circ \varphi_q$ just for elements $a$ generating $\C_q[x,y,z]$, which are $x$, $y$ and $z$. We have
\begin{align*}
  y^bx&=qy^{b-1}(q^{-1}yx+q^{-2}z^2-q^{-2}z^2)=q(q^{-1}yx+q^{-2}z^2)y^{b-1}-q^{-1}y^{b-1}z^2 \\
  &=yxy^{b-1}+q^{-1}z^2y^{b-1}-q^{-1}y^{b-1}z^2 \\
  &= (xy+q(q^2-q^{-2})z^2)y^{b-1}+q^{-1}z^2y^{b-1}-q^{-4b+3}z^2y^{b-1}\\
  &=xy^b+q^3z^2y^{b-1}-q^{-4b+3}z^2y^{b-1} =xy^b+q^{-2b+3}(q^{2b}-q^{-2b})z^2y^{b-1}
\intertext{and}
  yx^a&= q(q^{-1}yx+q^{-2}z^2-q^{-2}z^2)x^{a-1}=qx^{a-1}(q^{-1}yx+q^{-2}z^2)-q^{-1}z^2x^{a-1} \\
  &=x^{a-1}yx+q^{-1}x^{a-1}z^2-q^{-1}z^2x^{a-1} \\
  &=x^{a-1}(xy+q(q^2-q^{-2})z^2)+q^{-1}x^{a-1}z^2-q^{-4a+3}x^{a-1}z^2 \\
  &=x^ay+q^3x^{a-1}z^2-q^{-4a+3}x^{a-1}z^2= x^ay+q^{-2a+3}(q^{2a}-q^{-2a})x^{a-1}z^2
\end{align*}
for all $a,b \in \N$. However, the result holds for all $a,b \in \N_0$. Hence, we obtain
\begin{align*}
  \rho_x(x^az^cy^b)&=x^az^cy^bx=x^az^c(xy^b+q^{-2b+3}(q^{2b}-q^{-2b})z^2y^{b-1})\\
  &=q^{-2c}x^{a+1}z^cy^b+q^{-2b+3}(q^{2b}-q^{-2b})x^az^{c+2}y^{b-1}
\intertext{and}
  \lambda_y(x^az^cy^b)&=yx^az^cy^b=(x^ay+q^{-2a+3}(q^{2a}-q^{-2a})x^{a-1}z^2)z^cy^b \\
  &=q^{-2c}x^az^cy^{b+1}+q^{-2a+3}(q^{2a}-q^{-2a})x^{a-1}z^{c+2}y^b
\end{align*}
for all $a,b,c \in \N_0$, which implies
\begin{align*}
  \varphi_q^{-1} \circ (\lambda_x - \rho_x) \circ \varphi_q &= -q(q^2-q^{-2}) z^2\gamma_{q,y}^{-2}\partial_{q^2,y}+x-x\gamma_{q,z}^{-2} \\
  &=q(q-q^{-1})z\gamma_{q,x}^{-2}\big(-[2]_qz\gamma_{q,x}^2 \gamma_{q,y}^{-2} \partial_{q^2,y}+x\gamma_{q,x}^2\gamma_{q,z}^{-1}\partial_{q,z}\big) \\
  &=q(q-q^{-1})z\gamma_{q,x}^{-2}\sigma_q^*(F)
\intertext{and}
  \varphi_q^{-1} \circ (\lambda_y - \rho_y) \circ \varphi_q&= q(q^2-q^{-2}) z^2 \gamma_{q,x}^{-2}\partial_{q^2,x}-y+y\gamma_{q,z}^{-2} \\
  &=q(q-q^{-1})z\gamma_{q,y}^{-2}\big([2]_qz\gamma_{q,y}^2\gamma_{q,x}^{-2} \partial_{q^2,x}-y\gamma_{q,y}^2 \gamma_{q,z}^{-1} \partial_{q,z}\big) \\
  &=q(q-q^{-1})z\gamma_{q,y}^{-2} \sigma_q^*(E).
\end{align*}
Moreover, we have
\begin{align*}
  (\lambda_z-\rho_z)(x^az^cy^b)&=zx^az^cy^b-x^az^cy^bz=(q^{-2a}-q^{-2b})x^az^{c+1}y^b
\end{align*}
for all $a,b,c \in \N_0$, which gives us
\begin{align*}
\varphi_q^{-1} \circ (\lambda_z - \rho_z) \circ \varphi_q&= \gamma_{q,x}^{-2} - \gamma_{q,y}^{-2} = \gamma_{q,y}^{-2}\big(\gamma_{q,x}^{-2}\gamma_{q,y}^2-1\big)=\gamma_{q,y}^{-2}(\sigma_q^*(K)-1).
\end{align*}
Therefore, we have $f \in Z(\C_q[x,y,z])$ if and only if
\begin{gather*}
  q(q-q^{-1})z\gamma_{q,y}^{-2} \sigma_q^*(E)\varphi_q^{-1}(f)=0, \qquad  q(q-q^{-1})z\gamma_{q,x}^{-2}\sigma_q^*(F)\varphi_q^{-1}(f)=0, \\
   \gamma_{q,y}^{-2}(\sigma_q^*(K)-1)\varphi_q^{-1}(f)=0,
\end{gather*}
which is equivalent to
\begin{align*}
  \sigma_q^*(E)\varphi_q^{-1}(f)=0, \qquad  (\sigma_q^*(K)-1)\varphi_q^{-1}(f)=0, \qquad \sigma_q^*(F)\varphi_q^{-1}(f)=0.
\end{align*}
Hence, we obtain ${\rm Z}(\C_q[x,y,z])=\C_q[x,y,z]^{U_q(\mfrak{sl}(2,\C))}$. This finishes the proof.}

\proposition{We have
\begin{align}
  \Delta_q (p_q^{s+1})=[s+1]_{q^2}[s+{\textstyle {3 \over 2}}]_{q^2}p_q^s
\end{align}
for all $s \in \N_0$. Therefore, the polynomial $b_q(s)=[s+1]_{q^2}[s+{\textstyle {3 \over 2}}]_{q^2}$ is a $q$-analogue of the Bernstein--Sato polynomial.}

\proof{By Lemma \ref{lem:commutator U_q(sl(2,C))}, Proposition \ref{prop:invariant operators and dual group} and Lemma \ref{lem:tranformation of P_q} we have
\begin{align*}
  \Delta_q(p_q^{s+1})=[\Delta_q,p_q^{s+1}](1)=[s+1]_{q^2}[s+{\textstyle {3 \over 2}}]_{q^2}p_q^s
\end{align*}
for all $s \in \N_0$.}
\vspace{-2mm}


\section{The Howe duality $(U_q(\mfrak{sl}(2,\C)),U_q(\mfrak{sl}(n,\C)))$}
\label{section-third}

In this section we present another example of a $q$-deformed Howe dual pair $(\mfrak{sl}(2,\C),\mfrak{sl}(n,\C))$.
The presence of the quantum Howe duality was implicitly supposed form the very beginning of quantum group theory. A general quantum Howe duality (in particular in type $A$) was considered in \cite{Lehrer-Zhang-Zhang2011}. However, we describe a different kind of quantum Howe duality of type $A$. This particular example shows that we may expect similar results for $q$-deformed analogues of all classical Howe dual pairs.


\subsection{The quantum group $U_q(\mfrak{sl}(n,\C))$ and its representations}
Let $P$ be a free $\mathbb{Z}$-lattice of rank $n \geq 2$ with the
canonical basis $\{\veps_1,\dots,\veps_n\}$, i.e.\ $P=\bigoplus_{i=1}^n \mathbb{Z}\veps_i$, endowed with a symmetric bilinear form $(\veps_i,\veps_j)=\delta_{ij}$. Then $\Delta=\{\veps_i-\veps_j;\, 1\leq i\neq j\leq n\}$ is a root system of $\mfrak{sl}(n,\C)$. A positive root system in $\Delta$ is $\Delta^+ = \{\veps_i - \veps_j;\, 1 \leq i < j \leq n\}$ and the set of simple roots is $\Pi=\{\alpha_1,\alpha_2,\dots,\alpha_{n-1}\}$ with $\alpha_i=\veps_i-\veps_{i+1}$ for $i=1,2\dots,n-1$. The fundamental weights are $\omega_i = \smash{\sum_{j=1}^i} \veps_j$ for $i=1,2,\dots,n-1$.

Let $q \in \C^\times$ satisfies $q \neq \pm 1$. The quantum group $U_q(\mfrak{sl}(n,\C))$ is a unital associative $\C$-algebra generated by $E_i$, $F_i$, $K_i$, $K_i^{-1}$ for $i=1,2,\dots,n-1$ subject to the relations
\begin{align}
\begin{gathered}
    K_i K_i^{-1}=1, \qquad K_i K_j  = K_j K_i, \qquad    K_i^{-1}K_i=1, \\
    K_iE_jK_i^{-1} = q^{(\alpha_i,\alpha_j)}E_j, \qquad [E_i,F_j]= \delta_{ij} {K_i-K_i^{-1} \over q-q^{-1}}, \qquad K_iF_jK_i^{-1}= q^{-(\alpha_i,\alpha_j)}F_j, \\
\end{gathered}
\end{align}
and the quantum Serre relations
\begin{align}
\begin{gathered}
    E_i^2E_j-(q+q^{-1})E_iE_jE_i+E_jE_i^2=0  \qquad  (|i-j|=1),\\
     F_i^2F_j-(q+q^{-1})F_iF_jF_i+F_jF_i^2=0  \qquad  (|i-j|=1),\\
E_iE_j=E_jE_i,\qquad  F_iF_j=F_jF_i  \qquad (|i-j|>1)
\end{gathered}
\end{align}
for $i,j = 1,2,\dots,n-1$. Moreover, there is a unique Hopf algebra structure on $U_q(\mfrak{sl}(n,\C))$ with the coproduct $\Delta \colon U_q(\mfrak{sl}(n,\C)) \rarr U_q(\mfrak{sl}(n,\C)) \otimes_\C U_q(\mfrak{sl}(n,\C))$, the counit $\veps \colon U_q(\mfrak{sl}(n,\C)) \rarr \C$ and the antipode $S \colon U_q(\mfrak{sl}(n,\C)) \rarr U_q(\mfrak{sl}(n,\C))$ given by
\begin{align} \label{eq:Hopf algebra structure U_q(sl(n,C))}
\begin{gathered}
 \Delta(E_i)=E_i \otimes K_i + 1 \otimes E_i,\quad \Delta(K_i)=K_i \otimes K_i,\quad \Delta(F_i)=F_i \otimes 1 + K^{-1}_i \otimes F_i, \\
 \veps(E_i)=0, \qquad \qquad \veps(K_i) =1, \qquad \qquad \veps(F_i)=0, \\
 S(E_i)=-E_iK^{-1}_i, \qquad S(K_i)=K^{-1}_i, \qquad S(F_i)=-K_iF_i
\end{gathered}
\end{align}
for all $i=1,2,\dots,n-1$.

\medskip

Let $V=\C^n$ and $V^*=\C^n{}^*$ be the standard and the dual standard $U_q(\mfrak{sl}(n,\C))$-modules, respectively, defined by
\begin{gather}
\begin{gathered}
  \rho_{q,V}(E_i)=E_{i,i+1}, \qquad \rho_{q,V}(F_i)=E_{i+1,i}, \\
  \rho_{q,V}(K_i)=\sum_{j=1}^{i-1} E_{j,j} + qE_{i,i}+q^{-1}E_{i+1,i+1} + \sum_{j=i+2}^n E_{j,j}
\end{gathered} \label{eq:standard rep}
\intertext{with respect to the canonical basis $y=(y_1,y_2,\dots,y_n)$ of $V \simeq V^{**}$ and by}
\begin{gathered}
  \rho_{q,V^*}(E_i)=-E_{i+1,i}, \qquad   \rho_{q,V^*}(F_i)=-E_{i,i+1}, \\
  \rho_{q,V^*}(K_i)= \sum_{j=1}^{i-1} E_{j,j} +q^{-1}E_{i,i}+qE_{i+1,i+1} + \sum_{j=i+2}^n E_{j,j},
\end{gathered} \label{eq:dual standard rep}
\end{gather}
with respect to the canonical basis $x=(x_1,x_2,\dots,x_n)$ of $V^*$, for all $i=1,2,\dots,n-1$.   Here $E_{i,j}$ denotes the $(n \times n)$-matrix having $1$ at the intersection of the $i$-th row and $j$-th column and $0$ elsewhere.

Since we have $V\otimes_\C V\simeq L(2\omega_1) \oplus L(\omega_2)$ and $\smash{V^* \otimes_\C V^*} \simeq L(2\omega_{n-1}) \oplus L(\omega_{n-2})$, where $L(\lambda)$ is the simple $U_q(\mfrak{sl}(n,\C))$-module with highest weight $q^\lambda$ for $\lambda \in \mfrak{h}^*$, we define
\begin{align}
  \C_q[V] \simeq S_q(V^*) = T(V^*)/I_{q,V^*} \qquad \text{and} \qquad \C_q[V^*] \simeq S_q(V) = T(V)/I_{q,V},
\end{align}
where $I_{q,V}$ is the two-sided ideal of the tensor algebra $T(V)$ generated by $L(\omega_2)=\langle qy_i \otimes y_j - y_j \otimes y_i;\, 1 \leq i < j \leq n \rangle$ and
$I_{q,V^*}$ is the two-sided ideal of the tensor algebra $T(V^*)$ generated by $L(\omega_{n-2}) = \langle x_i \otimes x_j - qx_j \otimes x_i;\, 1 \leq i < j \leq j \rangle$, which gives us
\begin{align}
  \C_q[V]& \simeq \C_q[x] = \C \langle x \rangle/(x_ix_j-qx_jx_i;\, 1 \leq i < j \leq n) \label{eq:quantum ring dual fundamental}
\intertext{and}
  \C_q[V^*]& \simeq \C_q[y] = \C\langle y \rangle/(qy_iy_j-y_jy_i;\, 1 \leq i < j \leq n). \label{eq:quantum ring fundamental}
\end{align}
Moreover, since the two-sided ideals $I_{q,V}$ and $I_{q,V^*}$ are $U_q(\mfrak{sl}(n,\C))$-submodules, we obtain that $\C_q[V^*]$ and $\C_q[V]$ are $U_q(\mfrak{sl}(n,\C))$-modules. Taking the limit $q \rarr 1$, we get $I_{q,V} \rarr I_V$ and $I_{q,V^*} \rarr I_{V^*}$, hence $\C_q[y] \rarr \C[y]$ and $\C_q[x] \rarr \C[x]$.

Since the $\C$-algebras $\C_q[x]$ and $\C_q[y]$ have basis $\{x_1^{a_1}x_2^{a_2}\dots x_n^{a_n};\, a_1,a_2,\dots,a_n \in \N_0\}$ and $\{y_1^{b_1}y_2^{b_2}\dots y_n^{b_n};\, b_1,b_2,\dots,b_n \in \N_0\}$ respectively, we can find  families of isomorphisms $\varphi_{q,V} \colon \C[x] \rarr \C_q[x]$ and $\varphi_{q,V^*} \colon \C[y] \rarr \C_q[y]$ of vector spaces such that $\varphi_{q,V} \rarr \id$ and $\varphi_{q,V^*} \rarr \id $ for $q \rarr 1$. Let us define $\varphi_{q,V} \colon \C[x] \rarr \C_q[x]$ and $\varphi_{q,V^*} \colon \C[y] \rarr \C_q[y]$ by
\begin{align}
  \varphi_{q,V}(x_1^{a_1}x_2^{a_2}\dots x_n^{a_n})=x_1^{a_1}x_2^{a_2}\dots x_n^{a_n} \qquad \text{and} \qquad  \varphi_{q,V^*}(y_1^{b_1}y_2^{b_2}\dots y_n^{b_n})=y_1^{b_1}y_2^{b_2}\dots y_n^{b_n} \label{eq:isomomorphism polynomials U_q(sl(n,C))}
\end{align}
for all $a_1,a_2,\dots,a_n \in \N_0$ and $b_1,b_2,\dots,b_n \in \N_0$. We shall denote
\begin{align}
  x^a=x_1^{a_1}x_2^{a_2}\dots x_n^{a_n}\qquad \text{and} \qquad y^b=y_1^{b_1}y_2^{b_2}\dots y_n^{b_n}
\end{align}
for $a=(a_1,a_2,\dots,a_n) \in \N_0^n$ and $b=(b_1,b_2,\dots,b_n) \in \N_0^n$. Furthermore, we denote by $1_k \in \N_0^n$ for $k=1,2,\dots,n$ the $n$-tuple having $1$ in the $k$-th coordinate and $0$ otherwise.
Then the corresponding $U_q(\mfrak{sl}(n,\C))$-module structures on $\C[x]$ and $\C[y]$ are given through the homomorphisms
\begin{align}
\sigma_{q,V} \colon U_q(\mfrak{sl}(n,\C)) \rarr \eus{A}^q_V \qquad \text{and} \qquad \sigma_{q,V^*} \colon U_q(\mfrak{sl}(n,\C)) \rarr \eus{A}^q_{V^*}
\end{align}
of associative $\C$-algebras, where $\eus{A}^q_V$ and $\eus{A}^q_{V^*}$ are the quantum Weyl algebras of the vector spaces $V$ and $V^*$, defined by
\begin{align}
  \sigma_{q,V}(a) = \varphi_{q,V}^{-1} \circ \rho_{q,V^*}(a) \circ \varphi_{q,V} \qquad \text{and} \qquad  \sigma_{q,V^*}(a) = \varphi_{q,V^*}^{-1} \circ \rho_{q,V}(a) \circ \varphi_{q,V^*}
\end{align}
for all $a \in U_q(\mfrak{sl}(n,\C))$.  Let us note that a different choice of the isomorphisms $\varphi_{q,V}$ and $\varphi_{q,V^*}$ leads only to a distinct realization of the $U_q(\mfrak{sl}(n,\C))$-modules $\C_q[V]$ and $\C_q[V^*]$, respectively.
\medskip

\proposition{\label{prop:action and dual action U_q(sl(n,C))}
The homomorphisms $\sigma_{q,V} \colon U_q(\mfrak{sl}(n,\C)) \rarr \eus{A}^q_V$ and $\sigma_{q,V^*} \colon U_q(\mfrak{sl}(n,\C)) \rarr \eus{A}^q_{V^*}$ of associative $\C$-algebras are given by
\begin{gather}\label{eq:action U_q(sl(n,C))}
\begin{gathered}
  \sigma_{q,V}(E_i)=-x_{i+1} \gamma_{q,x_i}^{-1}\gamma_{q,x_{i+1}} \partial_{q,x_i}, \qquad
  \sigma_{q,V}(F_i)=-x_i \gamma_{q,x_i} \gamma_{q,x_{i+1}}^{-1} \partial_{q,x_{i+1}}, \\
  \sigma_{q,V}(K_i)=\gamma_{q,x_i}^{-1} \gamma_{q,x_{i+1}}
\end{gathered}
\intertext{and} \label{eq:dual action U_q(sl(n,C))}
\begin{gathered}
  \sigma_{q,V^*}(E_i)=y_i\partial_{q,y_{i+1}}, \qquad \sigma_{q,V^*}(F_i)=y_{i+1} \partial_{q,y_i}, \\
  \sigma_{q,V^*}(K_i)=\gamma_{q,y_i} \gamma_{q,y_{i+1}}^{-1}
\end{gathered}
\end{gather}
for $i=1,2,\dots,n-1$.}

\proof{The proof is a straightforward computation. Using \eqref{eq:standard rep} and \eqref{eq:Hopf algebra structure U_q(sl(n,C))}, we may write
\begin{align*}
  \rho_{q,V^*}(K_i)(x_1^{a_1} \dots x_i^{a_i}x_{i+1}^{a_{i+1}} \dots x_n^{a_n})=q^{-a_i +a_{i+1}} x_1^{a_1} \dots x_i^{a_i}x_{i+1}^{a_{i+1}} \dots x_n^{a_n},
\end{align*}
\begin{align*}
  \rho_{q,V^*}(E_i)(x_1^{a_1} \dots x_i^{a_i}x_{i+1}^{a_{i+1}} \dots x_n^{a_n})&=-\sum_{k=1}^{a_i} q^{-(a_i-k)+a_{i+1}} x_1^{a_1} \dots x_i^{k-1}x_{i+1}x_i^{a_i-k}x_{i+1}^{a_{i+1}} \dots x_n^{a_n} \\
  &=-\sum_{k=1}^{a_i} q^{-2(a_i-k)+a_{i+1}} x_1^{a_1} \dots x_i^{a_i-1}x_{i+1}^{a_{i+1}+1} \dots x_n^{a_n} \\
  &=-[a_i]_q q^{-a_i+a_{i+1}+1} x_1^{a_1} \dots x_i^{a_i-1}x_{i+1}^{a_{i+1}+1} \dots x_n^{a_n}
\intertext{and}
  \rho_{q,V^*}(F_i)(x_1^{a_1} \dots x_i^{a_i}x_{i+1}^{a_{i+1}} \dots x_n^{a_n})&=- \sum_{k=1}^{a_{i+1}} q^{a_i-(k-1)} x_1^{a_1} \dots x_i^{a_i} x_{i+1}^{k-1} x_i x_{i+1}^{a_{i+1}-k} \dots x_n^{a_n} \\
  &=- \sum_{k=1}^{a_{i+1}} q^{a_i-2(k-1)} x_1^{a_1} \dots x_i^{a_i+1} x_{i+1}^{a_{i+1}-1} \dots x_n^{a_n} \\
  &=-[a_{i+1}]_q q^{a_i-a_{i+1}+1} x_1^{a_1} \dots x_i^{a_i+1} x_{i+1}^{a_{i+1}-1} \dots x_n^{a_n}
\end{align*}
for all $a_1,a_2,\dots,a_n \in \N_0$, which gives us \eqref{eq:action U_q(sl(n,C))}.
Similarly, using \eqref{eq:dual standard rep} and \eqref{eq:Hopf algebra structure U_q(sl(n,C))}, we obtain
\begin{align*}
  \rho_{q,V}(K_i)(y_1^{b_1} \dots y_i^{b_i}y_{i+1}^{b_{i+1}} \dots y_n^{b_n})=q^{b_i-b_{i+1}} y_1^{b_1} \dots y_i^{b_i}y_{i+1}^{b_{i+1}} \dots y_n^{b_n},
\end{align*}
\begin{align*}
  \rho_{q,V}(E_i)(y_1^{b_1} \dots y_i^{b_i}y_{i+1}^{b_{i+1}} \dots y_n^{b_n})&= \sum_{k=1}^{b_{i+1}} q^{-b_{i+1}+k} y_1^{b_1} \dots y_i^{b_i} y_{i+1}^{k-1} y_i y_{i+1}^{b_{i+1}-k} \dots y_n^{b_n} \\
  &= \sum_{k=1}^{b_{i+1}} q^{-b_{i+1}+2(k-1)+1} y_1^{b_1} \dots y_i^{b_i+1} y_{i+1}^{b_{i+1}-1} \dots y_n^{b_n} \\
  &=[b_{i+1}]_q  y_1^{b_1} \dots y_i^{b_i+1} y_{i+1}^{b_{i+1}-1} \dots y_n^{b_n}
\intertext{and}
  \rho_{q,V}(F_i)(y_1^{b_1} \dots y_i^{b_i}y_{i+1}^{b_{i+1}} \dots y_n^{b_n}) &=\sum_{k=1}^{b_i} q^{-(k-1)} y_1^{b_1} \dots y_i^{k-1}y_{i+1}y_i^{b_i-k}y_{i+1}^{b_{i+1}} \dots y_n^{b_n} \\
  &=\sum_{k=1}^{b_i} q^{b_i-2(k-1)-1} y_1^{b_1} \dots y_i^{b_i-1}y_{i+1}^{b_{i+1}+1} \dots y_n^{b_n} \\
  &=[b_i]_q y_1^{b_1} \dots y_i^{b_i-1}y_{i+1}^{b_{i+1}+1} \dots y_n^{b_n}
\end{align*}
for all $b_1,b_2,\dots,b_n \in \N_0$, which implies \eqref{eq:dual action U_q(sl(n,C))}. This finishes the proof.}

Let us consider the tensor product $\C_q[V] \otimes_\C \C_q[V^*]$, which has the natural structure of a $U_q(\mfrak{sl}(n,\C))$-module. Our aim is to decompose this tensor product as a representation of $U_q(\mfrak{sl}(n,\C))$. For that that reason, let us introduce a unital associative $\C$-algebra
\begin{align}
  \C_q[x,y]=\C\langle x,y\rangle/(x_ix_j-qx_jx_i, qy_iy_j-y_jy_i,qx_iy_j-y_jx_i,x_jy_i-qy_ix_j;\,1 \leq i < j \leq n). \label{eq:quantum ring full}
\end{align}
Since $\C_q[x] \simeq \C_q[V]$ and $\C_q[y] \simeq \C_q[V^*]$ are subalgebras of $\C_q[x,y]$, we identify naturally
\begin{align}
  \C_q[V] \otimes_\C \C_q[V^*] \simeq \C_q[x] \otimes_\C \C_q[y] \riso \C_q[x,y] \label{eq:identification}
\end{align}
using the multiplication in $\C_q[x,y]$. As a result of this identification and
the actions $\rho_{q,V^*}$ on $\C_q[x]$, $\rho_{q,V}$ on $\C_q[y]$, we have an action $\rho_q$ of $U_q(\mfrak{sl}(n,\C))$ on $\C_q[x,y]$ through the coproduct defined by \eqref{eq:Hopf algebra structure U_q(sl(n,C))}.

As $\C[x,y]$ and $\C_q[x,y]$ are isomorphic vector spaces, we can again find a family of isomorphisms $\varphi_q \colon \C[x,y] \rarr \C_q[x,y]$ such that $\varphi_q \rarr \id $ for $q \rarr 1$. Let us define $\varphi_q \colon \C[x,y] \rarr \C_q[x,y]$ by
\begin{align}
  \varphi_q(x_1^{a_1}x_2^{a_2}\dots x_n^{a_n}y_1^{b_1}y_2^{b_2}\dots y_n^{b_n})=x_1^{a_1}x_2^{a_2}\dots x_n^{a_n}y_1^{b_1}y_2^{b_2}\dots y_n^{b_n}
\end{align}
for all $a_1,a_2,\dots,a_n \in \N_0$ and $b_1,b_2,\dots,b_n \in \N_0$. Then the corresponding $U_q(\mfrak{sl}(n,\C))$-module structure on $\C[x,y]$ is given through the homomorphism
\begin{align}
  \sigma_q \colon U_q(\mfrak{sl}(n,\C)) \rarr \eus{A}^q_{V \oplus V^*}
\end{align}
of associative $\C$-algebras, where $\eus{A}^q_{V \oplus V^*}$ is the quantum Weyl algebra of the vector space $V \oplus V^*$, defined by
\begin{align}
  \sigma_q(a) = \varphi_q^{-1}\! \circ \rho_q(a) \circ \varphi_q \label{eq:full action U_q(sl(n,C))}
\end{align}
for all $a \in U_q(\mfrak{sl}(n,\C))$.
\medskip

\proposition{The homomorphism $\sigma_q \colon U_q(\mfrak{sl}(n,\C)) \rarr \eus{A}^q_{V \oplus V^*}$ of associative $\C$-algebras is given by
\begin{align}
  \begin{aligned}
    \sigma_q(E_i)&=-x_{i+1} \gamma_{q,x_i}^{-1}\gamma_{q,x_{i+1}} \partial_{q,x_i} \gamma_{q,y_i} \gamma_{q,y_{i+1}}^{-1} + y_i\partial_{q,y_{i+1}}, \\
    \sigma_q(K_i)&=\gamma_{q,x_i}^{-1}\gamma_{q,x_{i+1}} \gamma_{q,y_i} \gamma_{q,y_{i+1}}^{-1}, \\
    \sigma_q(F_i)&=-x_i \gamma_{q,x_i} \gamma_{q,x_{i+1}}^{-1} \partial_{q,x_{i+1}} + \gamma_{q,x_i}\gamma_{q,x_{i+1}}^{-1} y_{i+1} \partial_{q,y_i}
  \end{aligned}
\end{align}
for $i=1,2,\dots,n-1$.}

\proof{It is a straightforward consequence of Proposition \ref{prop:action and dual action U_q(sl(n,C))}, the identification \eqref{eq:identification} and the formula \eqref{eq:full action U_q(sl(n,C))}.}
\vspace{-2mm}


\subsection{The Fischer decomposition for $U_q(\mfrak{sl}(n,\C))$}

Let us introduce the differential operators
\begin{align}
\begin{gathered}
  P_q= \sum_{k=1}^n q^{k-1}x_ky_k \prod_{i=1}^{k-1} \gamma_{q,y_i} \prod_{i=k+1}^n \gamma_{q,x_i}^{-1}, \qquad Q_q= \sum_{k=1}^n q^{-n+k} \partial_{q,x_k} \partial_{q,y_k} \prod_{i=1}^{k-1} \gamma_{q,x_i} \prod_{i=k+1}^n \gamma_{q,y_i}^{-1}, \\
  E_q^x = \prod_{k=1}^n \gamma_{q,x_k}, \qquad E_q= \prod_{k=1}^n \gamma_{q,x_k} \gamma_{q,y_k}, \qquad E^y_q = \prod_{k=1}^n \gamma_{q,y_k}
\end{gathered}
\end{align}
in the quantum Weyl algebra $\eus{A}^q_{V \oplus V^*}$.
\medskip

\proposition{\label{prop:Howe dual group sl(n,C))}
The differential operators $P_q$, $Q_q$, \smash{$E_q^x$}, \smash{$E_q^y$} and $E_q$ are $U_q(\mfrak{sl}(n,\C))$-invariant, i.e.\ they commute with $\sigma_q(a)$ for all $a \in U_q(\mfrak{sl}(n,\C))$. Moreover, the mapping $\pi_q \colon U_q(\mfrak{sl}(2,\C)) \rarr \eus{A}^q_{V\oplus V^*}$ uniquely determined by
\begin{align}
  \pi_q(E)=P_q, \qquad \pi_q(K)=q^nE_q, \qquad \pi_q(F)=-Q_q
\end{align}
gives rise to a homomorphism of associative $\C$-algebras.}

\proof{We may write
\begin{align*}
  [\sigma_q(E_i),P_q]&=\sum_{k=1}^n [y_i\partial_{q,y_{i+1}} -x_{i+1} \gamma_{q,x_i}^{-1}\gamma_{q,x_{i+1}} \partial_{q,x_i} \gamma_{q,y_i} \gamma_{q,y_{i+1}}^{-1}, q^{k-1}x_ky_k \prod_{i=1}^{k-1} \gamma_{q,y_i} \prod_{i=k+1}^n \gamma_{q,x_i}^{-1}]\\
  &=[y_i\partial_{q,y_{i+1}}, q^ix_{i+1}y_{i+1} \prod_{j=1}^i \gamma_{q,y_j} \prod_{j=i+2}^n \gamma_{q,x_j}^{-1}]\\
  &\quad-[x_{i+1} \gamma_{q,x_i}^{-1}\gamma_{q,x_{i+1}} \partial_{q,x_i} \gamma_{q,y_i} \gamma_{q,y_{i+1}}^{-1}, q^{i-1}x_iy_i \prod_{j=1}^{i-1} \gamma_{q,y_j} \prod_{j=i+1}^n \gamma_{q,x_j}^{-1}] \\
  &=\prod_{j=1}^{i-1}  \gamma_{q,y_j} \prod_{j=i+2}^n \gamma_{q,x_j}^{-1} [y_i\partial_{q,y_{i+1}}, q^ix_{i+1}y_{i+1} \gamma_{q,y_i}]\\
  &\quad -\prod_{j=1}^{i-1}  \gamma_{q,y_j} \prod_{j=i+2}^n \gamma_{q,x_j}^{-1} [x_{i+1} \gamma_{q,x_i}^{-1}\gamma_{q,x_{i+1}} \partial_{q,x_i} \gamma_{q,y_i} \gamma_{q,y_{i+1}}^{-1}, q^{i-1}x_iy_i  \gamma_{q,x_{i+1}}^{-1}] \\
  &=q^i \prod_{j=1}^i  \gamma_{q,y_j} \prod_{j=i+2}^n \gamma_{q,x_j}^{-1} (x_{i+1}y_i \gamma_{q,y_{i+1}}^{-1}- x_{i+1}y_i \gamma_{q,y_{i+1}}^{-1})=0,\\
  [\sigma_q(K_i),P_q]&=\sum_{k=1}^n [\gamma_{q,x_i}^{-1}\gamma_{q,x_{i+1}} \gamma_{q,y_i} \gamma_{q,y_{i+1}}^{-1}, q^{k-1}x_ky_k \prod_{j=1}^{k-1} \gamma_{q,y_j} \prod_{j=k+1}^n \gamma_{q,x_j}^{-1}]=0, \\
  [\sigma_q(F_i),P_q]&=\sum_{k=1}^n [\gamma_{q,x_i}\gamma_{q,x_{i+1}}^{-1} y_{i+1} \partial_{q,y_i} -x_i \gamma_{q,x_i} \gamma_{q,x_{i+1}}^{-1} \partial_{q,x_{i+1}}, q^{k-1}x_ky_k \prod_{j=1}^{k-1} \gamma_{q,y_j} \prod_{j=k+1}^n \gamma_{q,x_j}^{-1}] \\
  &=[\gamma_{q,x_i}\gamma_{q,x_{i+1}}^{-1} y_{i+1} \partial_{q,y_i}, q^{i-1}x_iy_i \prod_{j=1}^{i-1} \gamma_{q,y_j} \prod_{j=i+1}^n \gamma_{q,x_j}^{-1}] \\
  &\quad - [x_i \gamma_{q,x_i} \gamma_{q,x_{i+1}}^{-1} \partial_{q,x_{i+1}}, q^ix_{i+1}y_{i+1} \prod_{j=1}^i \gamma_{q,y_j} \prod_{j=i+2}^n \gamma_{q,x_j}^{-1}]\\
  &=\prod_{j=1}^{i-1} \gamma_{q,y_j} \prod_{j=i+2}^n \gamma_{q,x_j}^{-1} [\gamma_{q,x_i}\gamma_{q,x_{i+1}}^{-1} y_{i+1} \partial_{q,y_i}, q^{i-1}x_iy_i \gamma_{q,x_{i+1}}^{-1}] \\
  &\quad - \prod_{j=1}^{i-1} \gamma_{q,y_j} \prod_{j=i+2}^n \gamma_{q,x_j}^{-1} [x_i \gamma_{q,x_i} \gamma_{q,x_{i+1}}^{-1} \partial_{q,x_{i+1}}, q^ix_{i+1}y_{i+1} \gamma_{q,y_i}]\\
  &=q^i \prod_{j=1}^i \gamma_{q,y_j} \prod_{j=i+2}^n \gamma_{q,x_j}^{-1} (x_i y_{i+1} \gamma_{q,x_i} \gamma_{q,x_{i+1}}^{-2}- x_i y_{i+1} \gamma_{q,x_i} \gamma_{q,x_{i+1}}^{-2})=0
\end{align*}
and
\begin{align*}
  [\sigma_q(E_i),Q_q]&=\sum_{k=1}^n [y_i\partial_{q,y_{i+1}} -x_{i+1} \gamma_{q,x_i}^{-1}\gamma_{q,x_{i+1}} \partial_{q,x_i} \gamma_{q,y_i} \gamma_{q,y_{i+1}}^{-1}, q^{k-n} \partial_{q,x_k} \partial_{q,y_k} \prod_{j=1}^{k-1} \gamma_{q,x_j} \prod_{j=k+1}^n \gamma_{q,y_j}^{-1}]\\
  &=[y_i\partial_{q,y_{i+1}}, q^{i-n} \partial_{q,x_i} \partial_{q,y_i} \prod_{j=1}^{i-1} \gamma_{q,x_j} \prod_{j=i+1}^n \gamma_{q,y_j}^{-1}]\\
  &\quad -[x_{i+1} \gamma_{q,x_i}^{-1}\gamma_{q,x_{i+1}} \partial_{q,x_i} \gamma_{q,y_i} \gamma_{q,y_{i+1}}^{-1}, q^{i-n+1} \partial_{q,x_{i+1}} \partial_{q,y_{i+1}} \prod_{j=1}^i \gamma_{q,x_j} \prod_{j=i+2}^n \gamma_{q,y_j}^{-1}]\\
  &=\prod_{j=1}^{i-1} \gamma_{q,x_j} \prod_{j=i+2}^n \gamma_{q,y_j}^{-1} [y_i\partial_{q,y_{i+1}}, q^{i-n} \partial_{q,x_i} \partial_{q,y_i} \gamma_{q,y_{i+1}}^{-1}]\\
  &\quad - \prod_{j=1}^{i-1} \gamma_{q,x_j} \prod_{j=i+2}^n \gamma_{q,y_j}^{-1} [x_{i+1} \gamma_{q,x_i}^{-1}\gamma_{q,x_{i+1}} \partial_{q,x_i} \gamma_{q,y_i} \gamma_{q,y_{i+1}}^{-1}, q^{i-n+1} \partial_{q,x_{i+1}} \partial_{q,y_{i+1}} \gamma_{q,x_i}]\\
  &= -q^{i-n+1} \prod_{j=1}^{i-1} \gamma_{q,x_j} \prod_{j=i+1}^n \gamma_{q,y_j}^{-1} (\partial_{q,x_i} \partial_{q,y_{i+1}} \gamma_{q,y_i} - \partial_{q,x_i} \partial_{q,y_{i+1}} \gamma_{q,y_i})=0, \\
  [\sigma_q(K_i),Q_q]&=\sum_{k=1}^n [\gamma_{q,x_i}^{-1}\gamma_{q,x_{i+1}} \gamma_{q,y_i} \gamma_{q,y_{i+1}}^{-1}, q^{k-n} \partial_{q,x_k} \partial_{q,y_k} \prod_{j=1}^{k-1} \gamma_{q,x_j} \prod_{j=k+1}^n \gamma_{q,y_j}^{-1}]=0, \\
  [\sigma_q(F_i),Q_q]&=\sum_{k=1}^n [\gamma_{q,x_i}\gamma_{q,x_{i+1}}^{-1} y_{i+1} \partial_{q,y_i} -x_i \gamma_{q,x_i} \gamma_{q,x_{i+1}}^{-1} \partial_{q,x_{i+1}}, q^{k-n} \partial_{q,x_k} \partial_{q,y_k} \prod_{j=1}^{k-1} \gamma_{q,x_j} \prod_{j=k+1}^n \gamma_{q,y_j}^{-1}] \\
  &=[\gamma_{q,x_i}\gamma_{q,x_{i+1}}^{-1} y_{i+1} \partial_{q,y_i}, q^{i-n+1} \partial_{q,x_{i+1}} \partial_{q,y_{i+1}} \prod_{j=1}^i \gamma_{q,x_j} \prod_{j=i+2}^n \gamma_{q,y_j}^{-1}] \\
  &\quad -[x_i \gamma_{q,x_i} \gamma_{q,x_{i+1}}^{-1} \partial_{q,x_{i+1}}, q^{i-n} \partial_{q,x_i} \partial_{q,y_i} \prod_{j=1}^{i-1} \gamma_{q,x_j} \prod_{j=i+1}^n \gamma_{q,y_j}^{-1}] \\
  &=\prod_{j=1}^{i-1} \gamma_{q,x_j} \prod_{j=i+2}^n \gamma_{q,y_j}^{-1} [\gamma_{q,x_i}\gamma_{q,x_{i+1}}^{-1} y_{i+1} \partial_{q,y_i}, q^{i-n+1} \partial_{q,x_{i+1}} \partial_{q,y_{i+1}} \gamma_{q,x_i}] \\
  &\quad - \prod_{j=1}^{i-1} \gamma_{q,x_j} \prod_{j=i+2}^n \gamma_{q,y_j}^{-1} [x_i \gamma_{q,x_i} \gamma_{q,x_{i+1}}^{-1} \partial_{q,x_{i+1}}, q^{i-n} \partial_{q,x_i} \partial_{q,y_i} \gamma_{q,y_{i+1}}^{-1}] \\
  &= -q^{i-n+1}\prod_{j=1}^{i-1} \gamma_{q,x_j} \prod_{j=i+1}^n \gamma_{q,y_j}^{-1}  (\partial_{q,x_{i+1}}\partial_{q,y_i}\gamma_{q,x_i}^2 \gamma_{q,x_{i+1}}^{-1} - \partial_{q,x_{i+1}}\partial_{q,y_i}\gamma_{q,x_i}^2 \gamma_{q,x_{i+1}}^{-1}) =0
\end{align*}
for $i=1,2,\dots,n$. Furthermore, we immediately obtain
\begin{gather*}
  [\sigma_q(E_i),E^x_q]=0, \qquad [\sigma_q(K_i),E^x_q]=0, \qquad [\sigma_q(F_i),E^x_q]=0, \\
    [\sigma_q(E_i),E^y_q]=0, \qquad [\sigma_q(K_i),E^y_q]=0, \qquad [\sigma_q(F_i),E^y_q]=0.
\end{gather*}
Therefore, the differential operators $P_q$, $Q_q$, $E^x_q$, $E^y_q$ and $E_q$ are $U_q(\mfrak{sl}(n,\C))$-invariant, since we have $E_q=E^x_q E^y_q$. Moreover, we have
\begin{align*}
  \pi_q(K)\pi_q(E) &=q^n \prod_{k=1}^n \gamma_{q,x_k} \gamma_{q,y_k} \sum_{k=1}^n q^{k-1}x_ky_k \prod_{i=1}^{k-1} \gamma_{q,y_i} \prod_{i=k+1}^n \gamma_{q,x_i}^{-1}\\
  &=q^2\sum_{k=1}^n q^{k-1}x_ky_k \prod_{i=1}^{k-1} \gamma_{q,y_i} \prod_{i=k+1}^n \gamma_{q,x_i}^{-1} q^n \prod_{k=1}^n \gamma_{q,x_k} \gamma_{q,y_k} =q^2\pi_q(E)\pi_q(K), \\
 \pi_q(K)\pi_q(F) &=-q^n \prod_{k=1}^n \gamma_{q,x_k} \gamma_{q,y_k}\sum_{k=1}^n q^{-n+k} \partial_{q,x_k} \partial_{q,y_k} \prod_{i=1}^{k-1} \gamma_{q,x_i} \prod_{i=k+1}^n \gamma_{q,y_i}^{-1}\\
  &=-q^{-2} \sum_{k=1}^n q^{-n+k} \partial_{q,x_k} \partial_{q,y_k} \prod_{i=1}^{k-1} \gamma_{q,x_i} \prod_{i=k+1}^n \gamma_{q,y_i}^{-1} q^n \prod_{k=1}^n \gamma_{q,x_k} \gamma_{q,y_k}=q^{-2}\pi_q(F)\pi_q(K)
\end{align*}
and
\begin{align*}
  [\pi_q(E), \pi_q(F)]&=\sum_{k=1}^n\sum_{\ell=1}^n\Big[q^{-n+k} \partial_{q,x_k} \partial_{q,y_k} \prod_{i=1}^{k-1} \gamma_{q,x_i} \prod_{i=k+1}^n \gamma_{q,y_i}^{-1}, q^{\ell-1}x_\ell y_\ell \prod_{i=1}^{\ell-1} \gamma_{q,y_i} \prod_{i=\ell+1}^n \gamma_{q,x_i}^{-1}\Big] \\
  &=\sum_{k=1}^n\Big[q^{-n+k} \partial_{q,x_k} \partial_{q,y_k} \prod_{i=1}^{k-1} \gamma_{q,x_i} \prod_{i=k+1}^n \gamma_{q,y_i}^{-1}, q^{k-1}x_k y_k \prod_{i=1}^{k-1} \gamma_{q,y_i} \prod_{i=k+1}^n \gamma_{q,x_i}^{-1}\Big] \\
  &=\sum_{k=1}^nq^{-n+2k-1}\prod_{i=1}^{k-1} \gamma_{q,x_i} \gamma_{q,y_i} \prod_{i=k+1}^n \gamma_{q,x_i}^{-1} \gamma_{q,y_i}^{-1}[\partial_{q,x_k} \partial_{q,y_k},x_k y_k] \\
  &=\sum_{k=1}^nq^{-n+2k-1}\prod_{i=1}^{k-1} \gamma_{q,x_i} \gamma_{q,y_i} \prod_{i=k+1}^n \gamma_{q,x_i}^{-1} \gamma_{q,y_i}^{-1}{q \gamma_{q,x_k}\gamma_{q,y_k} - q^{-1}\gamma_{q,x_k}^{-1}\gamma_{q,y_k}^{-1} \over q-q^{-1}} \\
  &= {q^n \prod_{i=1}^n \gamma_{q,x_i} \gamma_{q,y_i} - q^{-n} \prod_{i=1}^n \gamma_{q,x_i}^{-1} \gamma_{q,y_i}^{-1} \over q-q^{-1}}= {\pi_q(K)-\pi_q(K)^{-1} \over q-q^{-1}},
\end{align*}
which gives rise to a homomorphism $\pi_q\colon U_q(\mfrak{sl}(2,\C)) \rarr \eus{A}^q_V$ of associative $\C$-algebras.}

\proposition{\label{prop:orthogonal decompostion U_q(sl(n,C))}
We have a decompostion
\begin{align}
  \C[x,y]= \bigoplus_{j=0}^\infty \bigoplus_{a,b=0}^\infty P_q^j H_{q,a,b},
\end{align}
where
\begin{align}
  H_{q,a,b}=\{f \in \C[x,y];\, Q_qf=0,\, E^x_qf=q^af,\, E^y_qf=q^bf\}
\end{align}
is the simple finite-dimensional highest weight $U_q(\mfrak{sl}(n,\C))$-module with highest weight $q^{b \omega_1 + a\omega_{n-1}}$ generated by the vector $x_n^ay_1^b \in \C[x,y]$ with
\begin{align}
\dim H_{q,a,b} =  {(a+b+n-1)(n-1) \over (a+n-1)(b+n-1)} \binom{a+n-1}{n-1}\!\binom{b+n-1}{n-1}
\end{align}
for all $a,b  \in \N_0$.}

\proof{First of all, we prove a decomposition
\begin{align}
  \C[x,y]_{a,b} = H_{q,a,b} \oplus P_q(\C[x,y]_{a-1,b-1}) \label{eq:direct sum U_q(sl(n,C))}
\end{align}
for all $a,b \in \N_0$, where $\C[x,y]_{a,b}$ denotes the common eigenspace of \smash{$E_q^x$} and \smash{$E_q^y$} on $\C[x,y]$ with eigenvalues $q^a$ and $q^b$, respectively. Let us suppose that $f \in  H_{q,a,b} \cap P_q(\C[x,y]_{a-1,b-1})$ and let us take the maximal integer $s \in \N$ such that $f = P_q^sg$ with $g \neq 0$. Using \eqref{eq:commutator U_q(sl(2,C))} and Proposition \ref{prop:Howe dual group sl(n,C))} we have
\begin{align*}
  0 = Q_q(f) = Q_q (P_q^s g) = [Q_q,P_q^s](g)+P^s_qQ_q(g)=P_q^sQ_q(g)+[s]_q[n+a+b-s-1]_qP_q^{s-1}(g),
\end{align*}
where we used the fact that $\pi_q(E)=P_q$, $\pi_q(K)=q^nE_q$ and $\pi_q(F)=-Q_q$. As $P_q$ is an injective differential operator, we obtain
\begin{align*}
  0=P_qQ_q(g)+[s]_q[n+a+b-s-1]_qg.
\end{align*}
Since the coefficient $[s]_q[n+a+b-s-1]_q$ is nonzero, we may write $f=P_q^{s+1}h$ with $h \neq 0$, which is in a contradiction with the maximality of the integer $s$. Therefore, we have $H_{q,a,b} \cap P_q(\C[x,y]_{a-1,b-1}) = \{0\}$. Further, since the mapping $Q_q \colon \C[x,y]_{a,b} \rarr \C[x,y]_{a-1,b-1}$ is injective on $P_q(\C[x,y]_{a-1,b-1}) \simeq \C[x,y]_{a-1,b-1}$, it is also surjective for all $a,b \in \N_0$, which implies that $\C[x,y]_{a,b}$ is a direct sum of $H_{q,a,b}$ and $P_q(\C[x,y]_{a-1,b-1})$ for all $a,b\in \N_0$.

As a consequence of the decomposition \eqref{eq:direct sum U_q(sl(n,C))} we immediately obtain a decomposition
\begin{align}
  \C[x,y]_{a,b}= \bigoplus_{j=0}^{\min\{a,b\}} P_q^j H_{q,a-j,b-j}, \label{eq:full decomposition U_q(sl(n,C))}
\end{align}
for all $a,b \in \N_0$. By Proposition \ref{prop:Howe dual group sl(n,C))} the differential operators $P_q$, \smash{$E^x_q$}, \smash{$E_q^y$} and $Q_q$ are $U_q(\mfrak{sl}(n,\C))$-invariant, therefore the subspaces $P_q^j H_{q,a,b}$ are $U_q(\mfrak{sl}(n,\C))$-submodules and $P_q^j H_{q,a,b} \simeq H_{q,a,b}$ as $U_q(\mfrak{sl}(n,\C))$-modules for all $j,a,b \in \in \N_0$. Finally, from \eqref{eq:direct sum U_q(sl(n,C))} we have
\begin{align*}
  \dim H_{q,a,b} &= \dim \C[x,y]_{a,b} - \dim \C[x,y]_{a-1,b-1} = {\textstyle \binom{a+n-1}{n-1}\binom{b+n-1}{n-1}}- {\textstyle \binom{a+n-2}{n-1}\binom{b+n-2}{n-1}} \\
  &={\textstyle {(a+b+n-1)(n-1) \over (a+n-1)(b+n-1)} \binom{a+n-1}{n-1}\binom{b+n-1}{n-1}}
\end{align*}
for all $a,b \in \N_0$. Since $\sigma_q(E_i)x_n^ay_1^b=0$,  $\sigma_q(K_i)x_n^ay_1^b=q^{b\omega_1+a\omega_{n-1}}(K_i)x_n^ay_1^b$ for $i=1,2,\dots,n-1$ and  $Q_q(x_n^ay_1^b)=0$, we obtain that $H_{q,a,b}$ contains a simple finite-dimensional highest weight $U_q(\mfrak{sl}(n,\C))$-submodule isomorphic to $L_q(b\omega_1+a\omega_{n-1})$, which is however isomorphic to $H_{q,a,b}$, because it has the same dimension as $H_{q,a,b}$ for all $a,b \in \N_0$.}

\lemma{\label{lem:tranformation of P_q U_q(sl(n,C))}
Let $p_q = \sum_{k=1}^n q^{k-1}x_ky_k$. We have
\begin{align}
  \varphi_q \circ P_q \circ \varphi_q^{-1} = p_q \qquad \text{and} \qquad \varphi_q \circ Q_q \circ \varphi_q^{-1} = \Delta_q,
\end{align}
where
\begin{align}
  \Delta_q(x^ay^b)=\sum_{k=1}^n q^{-n+k} q^{\sum_{i=1}^{k-1} a_i - \sum_{i=k+1}^n b_i} [a_k]_q[b_k]_q x^{a-1_k}y^{b-1_k}
\end{align}
for all $a,b \in \N_0^n$.}

\proof{We may write
\begin{align*}
  (\varphi_q \circ P_q \circ \varphi_q^{-1})(x^ay^b) &=\varphi_q(P_q(x^ay^b))=\sum_{k=1}^n q^{k-1}q^{\sum_{i=1}^{k-1} b_i -\sum_{i=k+1}^n a_i} \varphi_q(x^{a+1_k}y^{b+1_k}) \\
  &=\sum_{k=1}^n q^{k-1}q^{\sum_{i=1}^{k-1} b_i -\sum_{i=k+1}^n a_i} x^{a+1_k}y^{b+1_k} \\
  &=\sum_{k=1}^n q^{k-1}q^{\sum_{i=1}^{k-1} a_i -\sum_{i=k+1}^n a_i} x_kx^ay_ky^b = \sum_{k=1}^n q^{k-1}x_ky_k x^ay^b
\intertext{and}
  (\varphi_q \circ Q_q \circ \varphi_q^{-1})(x^ay^b)&=\varphi_q(Q_q(x^ay^b))= \sum_{k=1}^n q^{-n+k} q^{\sum_{i=1}^{k-1} a_i - \sum_{i=k+1}^n b_i} [a_k]_q[b_k]_q \varphi_q(x^{a-1_k}y^{b-1_k}) \\
  &= \sum_{k=1}^n q^{-n+k} q^{\sum_{i=1}^{k-1} a_i - \sum_{i=k+1}^n b_i} [a_k]_q[b_k]_q x^{a-1_k}y^{b-1_k}
\end{align*}
for all $a,b \in \N_0^n$.}

From Proposition \ref{prop:Howe dual group sl(n,C))} follows that $\C[x,y]$ is a $U_q(\mfrak{sl}(2,\C))$-module. Hence, also $\C_q[x,y]$ has a $U_q(\mfrak{sl}(2,\C))$-module structure given through the homomorphism
\begin{align}
  \tau_q\colon U_q(\mfrak{sl}(2,\C)) \rarr \End \C_q[x,y]
\end{align}
of associative $\C$-algebras, which is defined by
\begin{align}
  \tau_q(a)= \varphi_q \circ \pi_q(a) \circ \varphi_q^{-1}
\end{align}
for $a \in U_q(\mfrak{sl}(2,\C))$. Moreover, we have the original $U_q(\mfrak{sl}(n,\C))$-module structure on $\C_q[x,y]$, which commutes with the action of $U_q(\mfrak{sl}(2,\C))$ as follows from Proposition \ref{prop:Howe dual group sl(n,C))}. Therefore, we may decompose $\C_q[x,y]$ with respect to the action of $U_q(\mfrak{sl}(2,\C)) \otimes_\C U_q(\mfrak{sl}(n,\C))$. The corresponding decomposition is one of the main result of this section.
\medskip

We denote by
\begin{align}
  \mcal{I}_q = \C[p_q]
\end{align}
the algebra of quantum invariant polynomials generated by $p_q$ and by
\begin{align}
  \mcal{H}_q = \{f \in \C_q[x,y];\, \Delta_q f=0\}
\end{align}
the vector space of quantum harmonic polynomials. Furthermore, we have the decomposition
\begin{align}
  \mcal{H}_q = \bigoplus_{a,b=0}^\infty \mcal{H}_{q,a,b},
\end{align}
where $\mcal{H}_{q,a,b}$ is the vector space of quantum harmonic polynomials with eigenvalues $q^a$ and $q^b$ due to $E_q^x$ and $E^y_q$, respectively.
\medskip

\theorem{\label{thm:Fischer decomposition U_q(sl(n,C))}
We have a decomposition
\begin{align}
  \C_q[V] \otimes_\C \C_q[V^*] \simeq \mcal{I}_q \otimes_\C \mcal{H}_q = \bigoplus_{a,b=0}^\infty \mcal{I}_q \otimes_\C \mcal{H}_{q,a,b}.
\end{align}
Moreover, we have
\begin{align}
  \mcal{I}_q \otimes_\C \mcal{H}_{q,a,b} \simeq M_q((n+a+b)\omega) \otimes_\C L_q(b\omega_1+a\omega_{n-1})
\end{align}
as $(U_q(\mfrak{sl}(2,\C)) \otimes_\C U_q(\mfrak{sl}(n,\C)))$-modules, where $M_q(\lambda\omega)$ for $\lambda \in \C$ is the Verma module with lowest weight $q^{\lambda\omega}$ for $U_q(\mfrak{sl}(2,\C))$ and $L_q(b\omega_1+a\omega_{n-1})$ for $a,b \in \N_0$ is the simple finite-dimensional highest weight module with highest weight $q^{b\omega_1+a\omega_{n-1}}$ for $U_q(\mfrak{sl}(n,\C))$. Furthermore, the vector $1 \otimes x_n^ay_1^b \in \mcal{I}_q \otimes_\C \mcal{H}_{q,a,b}$ for $a,b \in \N_0$ is the (lowest, highest) weight vector with (lowest, highest) weight $(q^{(n+a+b)\omega}, q^{b\omega_1+a\omega_{n-1}})$.}

\proof{From Proposition \ref{prop:orthogonal decompostion U_q(sl(n,C))} we have the decomposition
\begin{align*}
  \C[x,y] = \bigoplus_{j=0}^\infty \bigoplus_{a,b=0}^\infty P_q^j H_{q,a,b}.
\end{align*}
Using the isomorphism $\varphi_q \colon \C[x,y] \rarr \C_q[x,y]$ and Lemma \ref{lem:tranformation of P_q U_q(sl(n,C))}, we obtain
\begin{align*}
  \C_q[x,y]= \bigoplus_{j=0}^\infty \bigoplus_{a,b=0}^\infty\, p_q^j \varphi_q(H_{q,a,b}) = \bigoplus_{j=0}^\infty \bigoplus_{a,b=0}^\infty\, p_q^j\mcal{H}_{q,a,b} \simeq \bigoplus_{a,b=0}^\infty \mcal{I}_q \otimes_\C \mcal{H}_{q,a,b},
\end{align*}
since $\mcal{H}_{q,a,b} = \varphi_q(H_{q,a,b})$ for all $a,b\in \N_0$. Further, we have
\begin{align*}
  \sigma_q(E)(P_q^jh)=P_q^j\sigma_q(E)h,\qquad \sigma_q(K)(P_q^jh)=P_q^j\sigma_q(K)h, \qquad \sigma_q(F)(P_q^jh)&=P_q^j\sigma_q(F)h
\end{align*}
and
\begin{align*}
  \pi_q(E)(P_q^jh)&=P_q^{j+1}h, \\
  \pi_q(K)(P_q^jh)&= q^{n+a+b+2j} P_q^jh,\\
  \pi_q(F)(P_q^jh)&=-Q_q(P_q^jh)=[P_q^j,Q_q](h)-P_q^jQ_q(h)=-[j]_q[n+a+b+j-1]_q P_q^{j-1}h
\end{align*}
for all $j \in \N_0$ and $h \in H_{q,a,b}$, where we used Lemma \ref{lem:commutator U_q(sl(2,C))} and Proposition \ref{prop:Howe dual group sl(n,C))}, which implies that the vector subspace $\bigoplus_{j=0}^\infty P_q^j H_{q,a,b}$ is isomorphic to $M_q((n+a+b)\omega) \otimes_\C L_q(b\omega_1+a\omega_{n-1})$ as a $(U_q(\mfrak{sl}(2,\C)) \otimes_\C U_q(\mfrak{sl}(n,\C)))$-module, since Proposition \ref{prop:orthogonal decompostion U_q(sl(n,C))} gives us that $H_{q,a,b}$ is isomorphic to $L_q(b \omega_1 + a\omega_{n-1})$ as $U_q(\mfrak{sl}(n,\C))$-module for $a,b \in \N_0$.}

\theorem{\label{thm:invariant polynomials and center U_q(sl(n,C))}
The algebra of invariant quantum polynomials $(\C_q[V] \otimes_\C \C_q[V^*])^{U_q(\mfrak{sl}(n,\C))}$ is given by
\begin{align}
  (\C_q[V] \otimes_\C \C_q[V^*])^{U_q(\mfrak{sl}(n,\C))} \simeq \mcal{I}_q=\C[p_q]
\end{align}
and the center ${\rm Z}(\C_q[V] \otimes_\C \C_q[V^*])$ by
\begin{align}
  {\rm Z}(\C_q[V] \otimes_\C \C_q[V^*]) \simeq \C.
\end{align}
}

\proof{By Theorem \ref{thm:Fischer decomposition U_q(sl(n,C))} we have
\begin{align*}
  (\C_q[V] \otimes_\C \C_q[V^*])^{U_q(\mfrak{sl}(n,\C))} \simeq \mcal{I}_q \otimes_\C \mcal{H}_{q,0,0} \simeq \mcal{I}_q,
\end{align*}
since $\mcal{H}_{q,0,0}$ is the trivial representation of $U_q(\mfrak{sl}(n,\C))$.

Further, let us denote by $\lambda_a$ and $\rho_a$ the left and right multiplication by an element $a \in \C_q[x,y]$, respectively. Since $\C_q[x,y]$ is generated by the subset $S=\{x_1,x_2,\dots,x_n,y_1,y_2,\dots,y_n\}$, the condition $f \in {\rm Z}(\C_q[x,y])$ is equivalent to $(\lambda_a-\rho_a)f=0$ for all $a \in S$. Since we have
\begin{align*}
  (\lambda_{x_k}-\rho_{x_k})x^ay^b&=x_kx^ay^b-x^ay^bx_k=q^{-\sum_{i=1}^{k-1} a_i} \Big(1-q^{\sum_{i=1}^{k-1}(a_i-b_i)-\sum_{i=k+1}^n(a_i-b_i)} \!\Big) x^{a+1_k}y^b
\intertext{and}
  (\lambda_{y_k}-\rho_{y_k})x^ay^b&=y_kx^ay^b-x^ay^by_k=q^{\sum_{i=k+1}^n b_i}\Big(q^{\sum_{i=1}^{k-1}(a_i+b_i)-\sum_{i=k+1}^n(a_i+b_i)} -1\!\Big) x^ay^{b+1_k}
\end{align*}
for all $k=1,2,\dots,n$, we obtain that $x^ay^b \in {\rm Z}(\C_q[x,y])$ for $a,b \in \N_0^n$ if and only if
\begin{align*}
  \sum_{i=1}^{k-1} (a_i-b_i) - \sum_{i=k+1}^n (a_i-b_i) = 0 \qquad \text{and} \qquad \sum_{i=1}^{k-1} (a_i+b_i) - \sum_{i=k+1}^n (a_i+b_i)=0
\end{align*}
for all $k=1,2,\dots,n$. By taking the sum and difference of the previous equations, we get
\begin{align*}
   \sum_{i=1}^{k-1} a_i - \sum_{i=k+1}^n a_i = 0 \qquad \text{and} \qquad  \sum_{i=1}^{k-1} b_i - \sum_{i=k+1}^n b_i = 0
\end{align*}
for all $k=1,2,\dots,n$, which has only the trivial solution $a=0$ and $b=0$, since $a,b \in \N_0^n$. Now, let $f \in {\rm Z}(\C_q[x,y])$ be of the form
\begin{align*}
  f = \sum_{a,b \in \N_0^n} c_{a,b} x^a y^b
\end{align*}
with $c_{a,b} \in \C$ for $a,b \in \N_0^n$. Since from the previous considerations we have $(\lambda_{x_k}-\rho_{x_k})x^ay^b=\alpha_{k,a,b} x^{a+1_k}y^b$ and $(\lambda_{y_k}-\rho_{y_k})x^ay^b=\beta_{k,a,b} x^ay^{b+1_k}$ for $a,b \in \N_0^n$, $k=1,2,\dots,n$ and some $\alpha_{k,a,b}, \beta_{k,a,b}\in \C$, we may write
\begin{align*}
  (\lambda_{x_k}-\rho_{x_k})f&= \sum_{a,b \in \N_0^n} c_{a,b} (\lambda_{x_k}-\rho_{x_k})x^ay^b= \sum_{a,b \in \N_0^n} c_{a,b} \alpha_{k,a,b} x^{a+1_k}y^b, \\
  (\lambda_{y_k}-\rho_{y_k})f&= \sum_{a,b \in \N_0^n} c_{a,b} (\lambda_{y_k}-\rho_{y_k})x^ay^b= \sum_{a,b \in \N_0^n} c_{a,b} \beta_{k,a,b} x^ay^{b+1_k}
\end{align*}
for all $k=1,2,\dots,n$. As $\{x^ay^b;\, a,b \in \N_0^n\}$ is a basis of $\C_q[x,y]$, we obtain that $f \in {\rm Z}(\C_q[x,y])$ if and only if $x^ay^b \in {\rm Z}(\C_q[x,y])$ for all $a,b \in \N_0^n$ satisfying $c_{a,b} \neq 0$. However, we have that $x^ay^b \in {\rm Z}(\C_q[x,y])$ only for $a=b=0 \in \N_0^n$. This finishes the proof.}

\proposition{We have
\begin{align}
  \Delta_q (p_q^{s+1})=[s+1]_q[s+n]_qp_q^s
\end{align}
for all $s \in \N_0$. Therefore, the polynomial $b_q(s)=[s+1]_q[s+n]_q$ is a $q$-analogue of the Bernstein--Sato polynomial.}

\proof{By Lemma \ref{lem:commutator U_q(sl(2,C))}, Proposition \ref{prop:Howe dual group sl(n,C))} and Lemma \ref{lem:tranformation of P_q U_q(sl(n,C))} we have
\begin{align*}
  \Delta_q(p_q^{s+1})=[\Delta_q,p_q^{s+1}](1)=[s+1]_q[s+n]_qp_q^s
\end{align*}
for all $s \in \N_0$.}

\vspace{-2mm}


\section{Verma modules and the tensor product decomposition problem}
\label{section-verma}

In this section we describe an explicit realization of the Verma module $M_q(\lambda\omega)$ with $\lambda \in \C$ for the quantum group $U_q(\mfrak{sl}(2,\C))$. We use this realization to decompose the tensor product of two Verma modules as a representation of  $U_q(\mfrak{sl}(2,\C))$.
\medskip

We use the notation introduced in Section \ref{section-strange}. 

Let $U_q(\mfrak{n})$ and $U_q(\widebar{\mfrak{n}})$ be the $\C$-subalgebras of $U_q(\mfrak{sl}(2,\C))$ generated by $E$ and $F$, respectively. The $\C$-subalgebra generated by $K$ will be denoted by $U_q(\mfrak{h})$.
\medskip

\theorem{\label{thm:Verma module}
Let $\lambda \in \C$.
\begin{enum}
\item[1)] The mapping $\pi_\lambda \colon U_q(\mfrak{sl}(2,\C)) \rarr \eus{A}^q_{\widebar{\mfrak{n}}}$ uniquely determined by
\begin{align}
  \pi_\lambda(E)=q^\lambda x^2 \partial_{q,x}+[\lambda]_q x \gamma_{q,x}^{-1}, \qquad \pi_\lambda(K)=q^\lambda\gamma_{q,x}^2, \qquad \pi_\lambda(F)=-\partial_{q,x}
\end{align}
is a homomorphism of associative $\C$-algebras. Moreover, the $\eus{A}^q_{\widebar{\mfrak{n}}}$-module $\C[\partial_{q,x}]$ is isomorphic to the Verma module $M_q(\lambda\omega)$ as a representation of $U_q(\mfrak{sl}(2,\C))$.
\smallskip

\item[2)] The mapping $\hat{\pi}_\lambda \colon U_q(\mfrak{sl}(2,\C)) \rarr \eus{A}^q_{\widebar{\mfrak{n}}^*}$ uniquely determined by
\begin{align}
  \hat{\pi}_\lambda(E)=q^\lambda y \partial_{q,y}^2 -[\lambda]_q \gamma_{q,y} \partial_{q,y}, \qquad
  \hat{\pi}_\lambda(K)=q^\lambda \gamma_{q,y}^{-2}, \qquad
  \hat{\pi}_\lambda(F)=-y
\end{align}
is a homomorphism of associative $\C$-algebras. Moreover, the $\eus{A}^q_{\widebar{\mfrak{n}}^*}$-module $\C[y]$ is isomorphic to the Verma module $M_q(\lambda\omega)$ as a representation of $U_q(\mfrak{sl}(2,\C))$.
\end{enum}}

\proof{1) We may write
\begin{align*}
  \pi_\lambda(K)\pi_\lambda(E)&=q^\lambda \gamma_{q,x}^2(q^\lambda x^2 \partial_{q,x}+[\lambda]_q x \gamma_{q,x}^{-1})= q^2(q^\lambda x^2 \partial_{q,x}+[\lambda]_q x \gamma_{q,x}^{-1})q^\lambda \gamma_{q,x}^2=q^2\pi_\lambda(E)\pi_\lambda(K), \\
  \pi_\lambda(K)\pi_\lambda(F)&=-q^\lambda \gamma_{q,x}^2 \partial_{q,x} = -q^{-2}\partial_{q,x} q^\lambda \gamma_{q,x}^2 = q^{-2} \pi_\lambda(F) \pi_\lambda(K)
\end{align*}
and
\begin{align*}
  [\pi_\lambda(E),\pi_\lambda(F)]&=[\partial_{q,x}, q^\lambda x^2 \partial_{q,x}+[\lambda]_q x \gamma_{q,x}^{-1}]=  q^\lambda [\partial_{q,x},x^2 \partial_{q,x}]+ [\lambda]_q [\partial_{q,x},x \gamma_{q,x}^{-1}] \\
  &=q^\lambda {\gamma_{q,x}^2-\gamma_{q,x}^{-2} \over q-q^{-1}} + {q^\lambda - q^{-\lambda} \over q-q^{-1}}\, \gamma_{q,x}^{-2} = {q^\lambda \gamma_{q,x}^2 - q^{-\lambda} \gamma_{q,x}^{-2}  \over q-q^{-1}} \\
  &={\pi_\lambda(K) - \pi_\lambda(K)^{-1} \over q-q^{-1}},
\end{align*}
which gives us a homomorphism $\pi_\lambda\colon U_q(\mfrak{sl}(2,\C)) \rarr \eus{A}^q_{\widebar{\mfrak{n}}}$ of associative $\C$-algebras. As $\C[\partial_{q,x}]$ has a canonical structure of an $\eus{A}^q_{\widebar{\mfrak{n}}}$-module, we get a $U_q(\mfrak{sl}(2,\C))$-module structure on $\C[\partial_{q,x}]$ through the homomorphism $\pi_\lambda \colon U_q(\mfrak{sl}(2,\C)) \rarr \eus{A}^q_{\widebar{\mfrak{n}}}$ of associative $\C$-algebras. Moreover, $\C[\partial_{q,x}]$ is a free $U_q(\widebar{\mfrak{n}})$-module of rank one with a free generator $1 \in \C[\partial_{q,x}]$ and $\pi_\lambda(K)1=q^\lambda$, which implies that $\C[\partial_{q,x}]$ is isomorphic to the Verma module $M_q(\lambda\omega)$ as a representation of $U_q(\mfrak{sl}(2,\C))$.

2) Using the quantum Fourier transform $\mcal{F} \colon \eus{A}^q_{\widebar{\mfrak{n}}} \rarr \eus{A}^q_{\widebar{\mfrak{n}}^*}$ defined by
\begin{align*}
  \mcal{F}(x)=-\partial_{q,y}, \qquad \mcal{F}(\partial_{q,x})=y, \qquad \mcal{F}(\gamma_{q,x})=q^{-1}\gamma_{q,y}^{-1},
\end{align*}
we obtain a homomorphism $\hat{\pi}_\lambda \colon U_q(\mfrak{sl}(2,\C)) \rarr \eus{A}^q_{\widebar{\mfrak{n}}^*}$ of associative $\C$-algebras given through
\begin{align*}
  \hat{\pi}_\lambda = \mcal{F} \circ \pi_{\lambda+2}
\end{align*}
for $\lambda \in \C$. Therefore, we have
\begin{align*}
  \hat{\pi}_\lambda(E)&= \mcal{F}(q^{\lambda+2} x^2 \partial_{q,x}+[\lambda+2]_q x \gamma_{q,x}^{-1})= q^{\lambda+2} \partial_{q,y}^2 y- [\lambda+2]_q q \partial_{q,y} \gamma_{q,y} = q^\lambda y \partial_{q,y}^2 -[\lambda]_q \gamma_{q,y} \partial_{q,y}, \\
  \hat{\pi}_\lambda(K)&= \mcal{F}(q^{\lambda+2} \gamma_{q,x}^2)=q^\lambda \gamma_{q,y}^{-2}, \\
  \hat{\pi}_\lambda(F)&= \mcal{F}(-\partial_{q,x})=-y.
\end{align*}
Further, since $\C[y]$ has a canonical structure of an $\eus{A}^q_{\widebar{\mfrak{n}}^*}$-module, we obtain a $U_q(\mfrak{sl}(2,\C))$-module structure on $\C[y]$ through the homomorphism $\hat{\pi}_\lambda \colon U_q(\mfrak{sl}(2,\C)) \rarr \eus{A}^q_{\widebar{\mfrak{n}}^*}$ of associative $\C$-algebras. Moreover, $\C[y]$ is a free $U_q(\widebar{\mfrak{n}})$-module of rank one with a free generator $1 \in \C[y]$ and $\hat{\pi}_\lambda(K)1=q^\lambda$, which gives us that $\C[y]$ is isomorphic to the Verma module $M_q(\lambda\omega)$ as a representation of $U_q(\mfrak{sl}(2,\C))$.}

\theorem{\label{thm:Verma module decomposition}
Let $\lambda \in \C$. Then the Verma module $M_q(\lambda\omega)$ is simple if $q^{2\lambda} \notin \{q^0,q^2,\dots\}$. Moreover, if $q^{2\lambda} \in \{q^0,q^2,\dots\}$, then we have the short exact sequence
\begin{align}
    0 \rarr M_q((-\lambda-2)\omega) \rarr M_q(\lambda\omega) \rarr L_q(\lambda\omega) \rarr 0
\end{align}
of $U_q(\mfrak{sl}(2,\C))$-modules, where $L_q(\lambda\omega)$ is the simple finite-dimensional $U_q(\mfrak{sl}(2,\C))$-module with highest weight $q^{\lambda\omega}$.}

Let us consider two Verma modules $M_q(\lambda\omega)$ and $M_q(\mu\omega)$ with $\lambda,\mu \in \C$. Then the tensor product $M_q(\lambda\omega) \otimes_\C M_q(\mu\omega)$ is isomorphic to $\C[x,y]$ with the action of $U_q(\mfrak{sl}(2,\C))$ given by
\begin{align}
\begin{aligned}
  \hat{\pi}_{\lambda,\mu}(E)&=(q^\lambda x \partial_{q,x}^2-[\lambda]_q \gamma_{q,x} \partial_{q,x})q^\mu \gamma_{q,y}^{-2} + q^\mu y \partial_{q,y}^2 - [\mu]_q \gamma_{q,y} \partial_{q,y}\\
  \hat{\pi}_{\lambda,\mu}(K)&= q^{\lambda+\mu} \gamma_{q,x}^{-2} \gamma_{q,y}^{-2}, \\
  \hat{\pi}_{\lambda,\mu}(F)&= -x -q^{-\lambda}y\gamma_{q,x}^2,
\end{aligned}
\end{align}
where we used Theorem \ref{thm:Verma module} and the coproduct $\Delta \colon U_q(\mfrak{sl}(2,\C)) \rarr  U_q(\mfrak{sl}(2,\C)) \otimes_\C  U_q(\mfrak{sl}(2,\C))$ defined by \eqref{eq:Hopf algebra structure U_q(sl(2,C))} to identify $\C[x] \otimes_\C \C[y]$ with $\C[x,y]$.

To decompose the tensor product of $M_q(\lambda\omega)$ and $M_q(\mu\omega)$ as a representation of the quantum group $U_q(\mfrak{sl}(2,\C))$, we need to find singular vectors in $M_q(\lambda\omega) \otimes_\C M_q(\mu\omega)$. We define a $U_q(\mfrak{h})$-module
\begin{align}
  \Sing (M_q(\lambda\omega) \otimes_\C M_q(\mu\omega))=\{v \in M_q(\lambda\omega) \otimes_\C M_q(\mu\omega);\, Ev=0\}
\end{align}
and call it the vector space of singular vectors. Since $M_q(\lambda\omega) \otimes_\C M_q(\mu\omega)$ is a semisimple $U_q(\mfrak{h})$-module, we obtain that $\Sing (M_q(\lambda\omega) \otimes_\C M_q(\mu\omega))$ is also a semisimple $U_q(\mfrak{h})$-module. Let us denote by $(M_q(\lambda\omega) \otimes_\C M_q(\mu\omega))_{\lambda+\mu-2n}$ and $\Sing(M_q(\lambda\omega) \otimes_\C M_q(\mu\omega))_{\lambda+\mu-2n}$ weight spaces of $M_q(\lambda\omega) \otimes_\C M_q(\mu\omega)$ and $\Sing(M_q(\lambda\omega) \otimes_\C M_q(\mu\omega))$ with weight $q^{(\lambda+\mu-2n)\omega}$ for $n \in \N_0$, respectively. It is enough to find all singular weight vectors.

Let us assume that a weight vector
\begin{align}
  v_n = \sum_{k=0}^n a_k x^{n-k} y^k,
\end{align}
where $a_k \in \C$ for $k=0,1,\dots,n$, with weight $q^{(\lambda+\mu-2n)\omega}$ is a singular vector in $\C[x,y]$, which means that $\hat{\pi}_{\lambda,\mu}(E)v_n=0$.
Since we have
\begin{align*}
  (q^\lambda  x \partial_{q,x}^2-[\lambda]_q \gamma_{q,x}\partial_{q,x})x^k&=(q^\lambda [k]_q[k-1]_q- [\lambda]_q[k]_qq^{k-1})x^{k-1} \\
  &= -[k]_q([\lambda]_q q^{k-1}-q^\lambda[k-1]_q)x^{k-1} \\
  &=-[k]_q[\lambda-k+1]_q x^{k-1}
\end{align*}
for all $k \in \N_0$, we may write
\begin{align*}
  \hat{\pi}_{\lambda,\mu}(E)v_n&=-\sum_{k=0}^n a_k q^{\mu-2k}[n-k]_q[\lambda-n+k+1]_q x^{n-k-1}y^k -\sum_{k=0}^n a_k[k]_q[\mu-k+1]_q x^{n-k}y^{k-1} \\
  &=-\sum_{k=1}^n \big(a_{k-1} q^{\mu+2-2k}[n-k+1]_q[\lambda-n+k]_q + a_k[k]_q[\mu-k+1]_q\big) x^{n-k}y^{k-1}.
\end{align*}
Threfore, the condition $\hat{\pi}_{\lambda,\mu}(E)v_n=0$ gives us the recurrence relation
\begin{align}
  a_k[k]_q[\mu-k+1]_q+a_{k-1} q^{\mu+2-2k}[n-k+1]_q[\lambda-n+k]_q=0 \label{eq:sing vectors recurrence}
\end{align}
for all $k=1,2\dots,n$. Hence, for a fixed $n \in \N_0$, the dimension of $\Sing(M_q(\lambda\omega) \otimes_\C M_q(\mu\omega))_{\lambda+\mu-2n}$ depends on the highest weights $q^{\lambda\omega}$ and $q^{\mu\omega}$.
\medskip

If $q^{2\lambda} \in q^{2\N_0}=\{q^0,q^2,\dots \}$ for $\lambda \in \C$, then we denote by $\lambda_{\rm int} \in \N_0$ the nonnegative integer satisfying $q^{2\lambda}=q^{2\lambda_{\rm int}}$.
\medskip

\lemma{\label{lem:singular vectors}
Let $\lambda,\mu \in \C$.
\begin{enum}
  \item[1)] If either $q^{2\lambda} \notin q^{2\N_0}$ or $q^{2\mu} \notin q^{2\N_0}$, then
  \begin{align}
    \Sing(M_q(\lambda\omega) \otimes_\C M_q(\mu\omega))_{\lambda+\mu-2n}= \C\, v^{\lambda,\mu}_n
  \end{align}
  for all $n \in \N_0$, where
  \begin{align*}
  v_n^{\lambda,\mu}=\sum_{k=0}^n (-1)^k {\mu -k\brack n-k}_q {\lambda-n+k \brack k}_q  q^{\mu k - k(k-1)} x^{n-k}y^k
\end{align*}
for $\lambda, \mu \in \C$ and $n \in \N_0$.
\smallskip

  \item[2)] If $q^{2\lambda} \in q^{2\N_0}$ and $q^{2\mu} \in q^{2\N_0}$, then
   \begin{align}
     \Sing(M_q(\lambda\omega) \otimes_\C M_q(\mu\omega))_{\lambda+\mu-2n} =
     \begin{cases}
       \C\, v^{\lambda,\mu}_n & \text{for $\lambda_{\rm int} \geq  n$ or $\mu_{\rm int} \geq n$}, \\[3mm]
       \C\, v^{\lambda,\mu}_{n,\pm} &  \begin{gathered}\!
       \text{for $\lambda_{\rm int}, \mu_{\rm int} \leq n-1$},  \\ \lambda_{\rm int} + \mu_{\rm int} < n-1, \hfill
        \end{gathered} \\[3mm]
       \C\, v^{\lambda,\mu}_{n,+} \oplus \C\, v^{\lambda,\mu}_{n,-} &  \begin{gathered}\!
       \text{for $\lambda_{\rm int}, \mu_{\rm int} \leq n-1$},  \\ \!\lambda_{\rm int} + \mu_{\rm int} \geq n-1, \hfill
        \end{gathered}
     \end{cases}
   \end{align}
   where
\begin{align*}
  v_{n,\pm}^{\lambda,\mu}=\sum_{k=\mu_{\rm int}+1}^{n-1-\lambda_{\rm int}} (-1)^k {([\lambda-n+\mu_{\rm int}+2]_q)_{k-\mu_{\rm int}-1} \over ([\mu_{\rm int}+2]_q)_{k-\mu_{\rm int}-1}} {([n-k+1]_q)_{k-\mu_{\rm int}-1} \over ([\mu-k+1]_q)_{k-\mu_{\rm int}-1}}\,  q^{\mu k - k(k-1)} x^{n-k}y^k
\end{align*}
   and
\begin{align*}
  v_{n,+}^{\lambda,\mu}&=\sum_{k=0}^{n-1-\lambda_{\rm int}} (-1)^k {\lambda-n+k \brack k}_q {([n-k+1]_q)_k \over ([\mu-k+1]_q)_k}\, q^{\mu k - k(k-1)} x^{n-k}y^k, \\
  v_{n,-}^{\lambda,\mu}&=\sum_{k=\mu_{\rm int} +1}^n (-1)^k {\mu -k\brack n-k}_q {([k+1]_q)_{n-k} \over ([\lambda-n+k+1]_q)_{n-k}}\, q^{\mu k - k(k-1)} x^{n-k}y^k.
\end{align*}
\end{enum}}

\proof{From the previous considerations, we know that a weight vector $v_n = \sum_{k=0}^n a_k x^{n-k} y^k$ with $a_k \in \C$ for $k=0,1,\dots,n$ is singular, if and only if the coefficients $a_k$ for $k=0,1,\dots,n$ satisfy the recurrence relation \eqref{eq:sing vectors recurrence}.
\smallskip

1) If either $q^{2\lambda} \notin q^{2\N_0}$ or $q^{2\mu} \notin q^{2\N_0}$, then we have either $[\lambda-n+k]_q \neq 0$ or  $[\mu-k+1]_q \neq 0$ for all $k=1,2,\dots,n$ and the recurrence relation \eqref{eq:sing vectors recurrence} has a unique solution with the initial condition either $a_n \in \C$ or $a_0 \in \C$. Therefore, we get $\dim \Sing(M_q(\lambda\omega) \otimes_\C M_q(\mu\omega))_{\lambda+\mu-2n}=1$.
\smallskip

2) Let us assume that $q^{2\lambda} \in q^{2\N_0}$ and $q^{2\mu} \in q^{2\N_0}$.

i) If either $\lambda_{\rm int} \geq n$ or $\mu_{\rm int} \geq n$, then we have again either $[\lambda-n+k]_q \neq 0$ or  $[\mu-k+1]_q \neq 0$ for all $k=1,2,\dots,n$ and the recurrence relation \eqref{eq:sing vectors recurrence} has a unique solution with the initial condition either $a_n \in \C$ or $a_0 \in \C$. Hence, we obtain $\dim \Sing(M_q(\lambda\omega) \otimes_\C M_q(\mu\omega))_{\lambda+\mu-2n}=1$.

ii) If $\lambda_{\rm int}, \mu_{\rm int} \leq n-1$ and $\lambda_{\rm int} + \mu_{\rm int} < n-1$, then the recurrence relation \eqref{eq:sing vectors recurrence} for $k=\mu_{\rm int}+1$ gives
\begin{align*}
  a_{\mu_{\rm int}} q^{\mu-2\mu_{\rm int}}[n-\mu_{\rm int}]_q[\lambda-n+\mu_{\rm int}+1]_q=0.
\end{align*}
Since $\lambda_{\rm int} + \mu_{\rm int} < n-1$, we get $[\lambda-n+\mu_{\rm int}+1]_q \neq 0$, which implies $a_{\mu_{\rm int}}=0$. As $[\lambda-n+k]_q \neq 0$ for all $k=1,2,\dots,\mu_{\rm int}$, using the recurrence relation \eqref{eq:sing vectors recurrence} for $k=1,2,\dots,\mu_{\rm int}$, we obtain $a_k=0$ for $k=1,2,\dots,\mu_{\rm int}$. Furthermore, we have $[\mu-k+1]_q \neq 0$ for $k=\mu_{\rm int}+2,\mu_{\rm int}+3,\dots,n$, hence \eqref{eq:sing vectors recurrence} for $k=\mu_{\rm int}+2,\mu_{\rm int}+3,\dots,n$ gives us a unique solution with the initial condition $a_{\mu_{\rm int}+1} \in \C$. Moreover, as $[\lambda-n+k]_q=0$ for $k=n-\lambda_{\rm int} > \mu_{\rm int}+1$, we get $a_k=0$ for $k=n-\lambda_{\rm int},n-\lambda_{\rm int}+1,\dots,n$. Therefore, we have $\dim \Sing(M_q(\lambda\omega) \otimes_\C M_q(\mu\omega))_{\lambda+\mu-2n}=1$.

iii) If $\lambda_{\rm int}, \mu_{\rm int} \leq n-1$ and $\lambda_{\rm int} + \mu_{\rm int} \geq n-1$, then the recurrence relation \eqref{eq:sing vectors recurrence} for $k=\mu_{\rm int}+1$ gives us
\begin{align*}
  a_{\mu_{\rm int}} q^{\mu-2\mu_{\rm int}}[n-\mu_{\rm int}]_q[\lambda-n+\mu_{\rm int}+1]_q=0.
\end{align*}
Hence, we get $a_{\mu_{\rm int}}=0$ or $[\lambda-n+\mu_{\rm int}+1]_q=0$. Since, we have $[\lambda-n+k]_q \neq 0$ for $k=n-\lambda_{\rm int}+1,n-\lambda_{\rm int}+2,\dots,n$, using the recurrence relation \eqref{eq:sing vectors recurrence} for $k=n-\lambda_{\rm int}+1,n-\lambda_{\rm int}+2,\dots,\mu_{\rm int}$, we get $a_k=0$ for $k=n-\lambda_{\rm int},n-\lambda_{\rm int}+1,\dots,\mu_{\rm int}$. As \eqref{eq:sing vectors recurrence} is satisfied for $k=n-\lambda_{\rm int}$ and $[\lambda-n+k]_q \neq 0$ for $k=1,2,\dots,n-\lambda_{\rm int}-1$, the recurrence relation \eqref{eq:sing vectors recurrence} for $k=1,2,\dots,n-\lambda_{\rm int}-1$ gives us a unique solution with the initial condition $a_{n-\lambda_{\rm int}-1} \in \C$. On the other hand, as \eqref{eq:sing vectors recurrence} is also satisfied for $k=\mu_{\rm int}+1$ and $[\mu-k+1]_q \neq 0$ for $k=\mu_{\rm int}+2,\mu_{\rm int}+3,\dots,n$, the recurrence relation \eqref{eq:sing vectors recurrence} for $k=\mu_{\rm int}+2,\mu_{\rm int}+3,\dots,n$ gives us a unique solution with the initial condition $a_{\mu_{\rm int}+1} \in \C$. Hence, we obtain $\dim \Sing(M_q(\lambda\omega) \otimes_\C M_q(\mu\omega))_{\lambda+\mu-2n}=2$.
\smallskip

The last step is to verify that the given polynomials $v^{\lambda,\mu}_n$, $v^{\lambda,\mu}_{n,\pm}$, $v^{\lambda,\mu}_{n,+}$ and $v^{\lambda,\mu}_{n,-}$ are nonzero and satisfy the recurrence relation \eqref{eq:sing vectors recurrence}, which is a straightforward computation.}

Let $q^{2\lambda} \in q^{2\N_0}$ for $\lambda \in \C$. Then $\dim \Ext_{\mcal{O}_q}(M_q(-(\lambda+2)\omega),M_q(\lambda\omega)) =1$ and therefore there exists a nontrivial extension $P_q(\lambda\omega)$ uniquely determined by the nonsplit short exact sequence
\begin{align}
  0 \rarr M_q(\lambda\omega) \rarr P_q(\lambda\omega) \rarr M_q(-(\lambda+2)\omega) \rarr 0
\end{align}
of $U_q(\mfrak{sl}(2,\C))$-modules. The extension $P_q(\lambda\omega)$ can be described as follows. We define
\begin{align}
  P_q(\lambda\omega)=\langle v_{\lambda,n},  w_{-\lambda-2,n}; n \in \N_0 \rangle
\end{align}
with the action of $U_q(\mfrak{sl}(2,\C))$ given by
\begin{align}
\begin{aligned}
  F v_{\lambda,n} & =v_{\lambda,n+1}, &\quad Fw_{-\lambda-2,n} &= w_{-\lambda-2,n+1}\\
  K v_{\lambda,n} &= q^{\lambda-2n}v_{\lambda,n}, &\quad K w_{-\lambda-2,n} &= q^{-\lambda-2-2n}  w_{-\lambda-2,n} \\
  E v_{\lambda,n} &= [n]_q[\lambda-n+1]_q v_{\lambda,n-1}, &\quad E w_{-\lambda-2,n} &= v_{\lambda,\lambda+n} - [n]_q[\lambda+n+1]_q w_{-\lambda-2,n-1}
\end{aligned}
\end{align}
for all $n \in \N_0$.

Let us recall that the formal character of the Verma module $M_q(\lambda\omega)$ with highest weight $q^{\lambda\omega}$ for $\lambda \in \C$ is given by
\begin{align}
  \ch M_q(\lambda\omega) = \sum_{n=0}^\infty \dim M_q(\lambda\omega)_{\lambda-2n} w^{\lambda-2n} = {w^\lambda  \over 1 -w^{-2}},
\end{align}
where $w$ is a formal variable. Moreover, for the formal character of $P_q(\lambda\omega)$ we have
\begin{align}
  \ch P_q(\lambda\omega) = \ch M_q(\lambda\omega) + \ch M_q(-(\lambda+2)\omega)
\end{align}
provided $q^{2\lambda} \in q^{2\N_0}$.
\medskip

\lemma{\label{lem:character decomposition}
We have
\begin{align}
  \ch (M_q(\lambda\omega) \otimes_\C M_q(\mu\omega)) = \sum_{n=0}^\infty \ch M_q((\lambda+\mu-2n)\omega)
\end{align}
for all $\lambda, \mu \in \C$.}

\proof{We may write
\begin{align*}
  \sum_{n=0}^\infty \ch M_q((\lambda+\mu-2n)\omega)&= \sum_{n=0}^\infty {w^{\lambda+\mu-2n} \over 1-w^{-2}} = {w^\lambda \over 1-w^{-2}} \sum_{n=0}^\infty {w^{\mu-2n}} = {w^\lambda \over 1-w^{-2}} {w^\mu \over 1-w^{-2}} \\
  & = \ch M_q(\lambda\omega) \ch M_q(\mu\omega) = \ch(M_q(\lambda\omega) \otimes_\C M_q(\mu\omega)),
\end{align*}
which gives the required statement.}

The previous lemma gives us a suggestion for the tensor product decomposition of two Verma modules for $U_q(\mfrak{sl}(2,\C))$. The exact claim is the content of the following theorem.
\medskip

Let us introduce a central element of $U_q(\mfrak{sl}(2,\C))$, called the quantum Casimir element, by the formula
\begin{align}
  \Cas_q = FE + { qK+q^{-1}K^{-1} \over (q-q^{-1})^2}.
\end{align}
It is easy to verify that $\Cas_q$ is a central element of $U_q(\mfrak{sl}(2,\C))$.
\medskip

\proposition{The quantum Casimir element $\Cas_q$ acts on the Verma module $M_q(\lambda\omega)$ as $c_\lambda\id_{M_q(\lambda\omega)}$, where
\begin{align}
  c_\lambda= { q^{\lambda+1}+q^{-\lambda-1} \over (q-q^{-1})^2} .
\end{align}
Moreover, we have $c_\lambda = c_\mu$ for some $\lambda,\mu \in \C$ if and only if $q^\mu=q^\lambda$ or $q^\mu=q^{-\lambda-2}$.}

\proof{Since the Verma module $M_q(\lambda\omega)$ for $\lambda \in \C$ is generated by the highest weight vector $v_\lambda$, it is enough to compute the action of $\Cas_q$ on $v_\lambda$. We may write
\begin{align*}
  \Cas_q v_\lambda & =
    FE  v_\lambda + { qK+q^{-1}K^{-1} \over (q-q^{-1})^2}\, v_\lambda
    ={ q^{\lambda+1}+q^{-\lambda-1} \over (q-q^{-1})^2}\, v_\lambda.
\end{align*}
Further, the equality $c_\lambda=c_\mu$ implies $q^{\lambda+1}+q^{-\lambda-1}= q^{\mu+1}+q^{-\mu-1}$, which gives us
\begin{align*}
  q^{\lambda+1}+q^{-\lambda-1}- q^{\mu+1}-q^{-\mu-1} = (q^{\lambda+1}-q^{\mu+1})(1-q^{-\lambda-\mu-2})=0.
\end{align*}
Hence, we obtain either $q^\mu=q^\lambda$ or $q^\mu=q^{-\lambda-2}$.}

Let us denote by
\begin{align*}
  S_{\lambda,\mu}=
  \begin{cases}
    \emptyset & \text{if $q^{2(\lambda+\mu)} \notin q^{2\N_0}$}, \\[2mm]
   \big\{0,1,\dots,\lfloor {(\lambda+\mu)_{\rm int} \over 2} \rfloor\big\} & \text{if $q^{2(\lambda+\mu)} \in q^{2\N_0}$, $q^{2\lambda} \notin q^{2\N_0}$ or $q^{2\mu} \notin q^{2\N_0}$}, \\[2mm]
   \big\{\min\{\lambda_{\rm int},\mu_{\rm int}\}+1,\dots,\lfloor {\lambda_{\rm int}+\mu_{\rm int} \over 2} \rfloor\big\} & \text{if $q^{2\lambda} \in q^{2\N_0}$, $q^{2\mu} \in q^{2\N_0}$},
  \end{cases}
\end{align*}
\begin{align*}
  S^c_{\lambda,\mu}=
  \begin{cases}
    \emptyset & \text{if $q^{2(\lambda+\mu)} \notin q^{2\N_0}$}, \\[2mm]
   \big\{\lfloor {(\lambda+\mu)_{\rm int} +1 \over 2} \rfloor +1,\dots,(\lambda+\mu)_{\rm int}+1\big\} & \text{if $q^{2(\lambda+\mu)} \in q^{2\N_0}$, $q^{2\lambda} \notin q^{2\N_0}$ or $q^{2\mu} \notin q^{2\N_0}$}, \\[2mm]
   \big\{\lfloor {\lambda_{\rm int}+\mu_{\rm int} +1 \over 2} \rfloor +1,\dots, \max\{\lambda_{\rm int},\mu_{\rm int}\}\big\} & \text{if $q^{2\lambda} \in q^{2\N_0}$, $q^{2\mu} \in q^{2\N_0}$}
  \end{cases}
\end{align*}
and
\begin{align*}
  R_{\lambda,\mu} = \N_0 \setminus (S_{\lambda,\mu} \cup S^c_{\lambda,\mu})
\end{align*}
subsets of $\N_0$ for all $\lambda, \mu \in \C$. Furthermore, we denote by
\begin{align*}
  T_{\lambda,\mu}=\{0,1,\dots,\lfloor {\textstyle {\lambda_{\rm int}+\mu_{\rm int} \over 2}} \rfloor\}, \qquad T_{\lambda,\mu}^c=\big\{\lfloor {\textstyle {\lambda_{\rm int}+\mu_{\rm int} +1 \over 2}} \rfloor +1,\dots,\lambda_{\rm int}+\mu_{\rm int}+1\big\}
  \end{align*}
  and
  \begin{align*}
    U_{\lambda,\mu} = \N_0 \setminus (T_{\lambda,\mu} \cup T^c_{\lambda,\mu})
  \end{align*}
subsets of $\N_0$ for all $\lambda,\mu \in \C$ satisfying $q^{2\lambda} \in q^{2\N_0}$ and $q^{2\mu} \in q^{2\N_0}$.
\medskip

\theorem{\label{thm:tensor product Verma quantum}
We have
  \begin{align}
    M_q(\lambda\omega) \otimes_\C M_q(\mu\omega) \simeq \bigoplus_{n \in S_{\lambda,\mu}} P_q((\lambda+\mu-2n)\omega) \oplus \bigoplus_{n \in R_{\lambda,\mu}} M_q((\lambda+\mu-2n)\omega).
  \end{align}
for all $\lambda,\mu \in \C$.}

\proof{Let us assume that $v_n \in \Sing(M_q(\lambda\omega) \otimes_\C M_q(\mu\omega))_{\lambda+\mu-2n}$ is a singular vector with weight $q^{(\lambda+\mu-2n)\omega}$ for some $n\in \N_0$. Then $v_n$ generates the highest weight submodule of $M_q(\lambda\omega) \otimes_\C M_q(\mu\omega)$ with highest weight $q^{(\lambda+\mu-2n)\omega}$ isomorphic either to $M_q((\lambda+\mu-2n)\omega)$ or to $L_q((\lambda+\mu-2n)\omega)$ as follows from Theorem \ref{thm:Verma module decomposition}. For $q^{2(\lambda+\mu-2n)} \notin q^{2\N_0}$ we have $M_q((\lambda+\mu-2n)\omega) \simeq L_q((\lambda+\mu-2n)\omega)$. On the other hand, if $q^{2(\lambda+\mu-2n)} \in q^{2\N_0}$, then $L_q((\lambda+\mu-2n)\omega)$ is a finite-dimensional simple module containing the lowest weight vector with weight $q^{-(\lambda+\mu-2n)\omega}$, which is annihilated by $F$. However, the element $F$ acts freely on $M_q(\lambda\omega) \otimes_\C M_q(\mu\omega)$, hence the singular vector $v_n$ can not generate a finite-dimensional simple submodule isomorphic to $L_q((\lambda+\mu-2n)\omega)$. Therefore, any singular vector $v_n \in \Sing(M_q(\lambda\omega) \otimes_\C M_q(\mu\omega))_{\lambda+\mu-2n}$ for $n\in \N_0$ generates a submodule isomorphic to the Verma module $M_q((\lambda+\mu-2n)\omega)$.

If $v_n$ is a singular vector in $M_q(\lambda\omega) \otimes_\C M_q(\mu\omega)$ with weight $q^{(\lambda+\mu-2n)\omega}$ for $n\in \N_0$, then $\Cas_q v_n = c_{\lambda+\mu-2n} v_n$. Let us denote by $N$ the submodule of $M_q(\lambda\omega) \otimes_\C M_q(\mu\omega)$ generated by the vector space of singular vectors $\Sing(M_q(\lambda\omega) \otimes_\C M_q(\mu\omega))$. Then we have an eigenspace decomposition
\begin{align}
  N = \bigoplus_{\nu \in \C} N_\nu
\end{align}
of $N$ with respect to the quantum Casimir element $\Cas_q$, where $N_\nu$ is the eigenspace with eigenvalue $c_\nu$ for $\nu \in \C$. Since $\Cas_q$ is a central element of $U_q(\mfrak{sl}(2,\C))$, the eigenspaces $N_\nu$ are $U_q(\mfrak{sl}(2,\C))$-submodules for all $\nu \in \C$. Moreover, for all $\nu \in \C$ the submodule $N_\nu$ is generated by singular vectors of weight $q^{\nu\omega}$. Furthermore, we have a short exact sequence
\begin{align}
  0 \rarr N \rarr M_q(\lambda\omega) \otimes_\C M_q(\mu\omega) \rarr (M_q(\lambda\omega) \otimes_\C M_q(\mu\omega)) /N \rarr 0  \label{eq:short exact sequence}
\end{align}
of $U_q(\mfrak{sl}(2,\C))$-modules.
\smallskip

1) If $q^{2(\lambda+\mu)} \notin q^{2\N_0}$, then the dimension of $\Sing(M_q(\lambda\omega) \otimes_\C M_q(\mu\omega))_{\lambda+\mu-2n}$ is $1$ for all $n \in \N_0$ as follows from Lemma \ref{lem:singular vectors} (i). Since $c_{\lambda+\mu-2m} \neq c_{\lambda+\mu-2n}$ for $m >  n$, we obtain
\begin{align*}
  N = \bigoplus_{n =0}^\infty  N_{\lambda+\mu-2n}.
\end{align*}
As $N_{\lambda+\mu-2n}$ is generated by $v^{\lambda,\mu}_n$ for $n \in \N_0$, we get $N_{\lambda+\mu-2n} \simeq M_q((\lambda+\mu-2n)\omega)$ for all $n \in \N_0$. Since $\ch(M_q(\lambda\omega) \otimes_\C M_q(\mu\omega)) = \ch N$, the short exact sequence \eqref{eq:short exact sequence} gives us $M_q(\lambda\omega) \otimes_\C M_q(\mu\omega) =N$.

2) If $q^{2(\lambda+\mu)} \in q^{2\N_0}$ and either $q^{2\lambda} \notin q^{2\N_0}$ or $q^{2\mu} \notin q^{2\N_0}$, then the dimension of $\Sing(M_q(\lambda\omega) \otimes_\C M_q(\mu\omega))_{\lambda+\mu-2n}$ is $1$ for all $n \in \N_0$ as follows from Lemma \ref{lem:singular vectors} (i). Since $c_{\lambda+\mu-2m}=c_{\lambda+\mu-2n}$ for $m > n$ if and only if $m=(\lambda+\mu)_{\rm int}-n+1$ for $n=0,1,\dots,\smash{\lfloor {(\lambda+\mu)_{\rm int} \over 2} \rfloor}$, we obtain
\begin{align*}
  N = \bigoplus_{n \in S_{\lambda,\mu}} N_{\lambda+\mu-2n} \oplus \bigoplus_{n \in R_{\lambda,\mu}} N_{\lambda+\mu-2n},
\end{align*}
where $N_{\lambda+\mu-2n}$ is generated by $v^{\lambda,\mu}_n$ and $v^{\lambda,\mu}_m$ with $m=(\lambda+\mu)_{\rm int}-n+1$ for $n \in S_{\lambda,\mu}$ and by $v^{\lambda,\mu}_n$ for $n \in R_{\lambda,\mu}$. As $N_{\lambda+\mu-2n}$ is generated by $v^{\lambda,\mu}_n$ for $n \in R_{\lambda,\mu}$, we get $N_{\lambda+\mu-2n} \simeq M_q((\lambda+\mu-2n)\omega)$ for all $n \in R_{\lambda,\mu}$. For $n \in S_{\lambda,\mu}$ the submodule of $N_{\lambda+\mu-2n}$ generated by $v^{\lambda,\mu}_n$ is isomorphic to $M_q((\lambda+\mu-2n)\omega)$, which is a reducible $U_q(\mfrak{sl}(2,\C))$-module containing a singular vector of weight $q^{-(\lambda+\mu-2n+2)\omega}=q^{(\lambda+\mu-2m)\omega}$ with $m=(\lambda+\mu)_{\rm int}-n+1$ as follows from Theorem \ref{thm:Verma module decomposition}. Moreover, this singular vector is a multiple of $v^{\lambda,\mu}_m$, since $N_{\lambda+\mu-2n}$ contains all singular vectors of weight $q^{(\lambda+\mu-2m)\omega}$. Therefore, we get $N_{\lambda+\mu-2n} \simeq M_q((\lambda+\mu-2n)\omega)$ for all $n \in S_{\lambda,\mu}$. Furthermore, we have
\begin{align*}
  \ch ((M_q(\lambda\omega) \otimes_\C M_q(\mu\omega))/N) = \sum_{n \in S^c_{\lambda,\mu}} \ch M_q((\lambda+\mu-2n)\omega),
\end{align*}
which gives us $\dim (M_q(\lambda\omega) \otimes_\C M_q(\mu\omega))/N)_{\lambda+\mu-2n}=n-\lfloor {(\lambda+\mu)_{\rm int}+1 \over 2} \rfloor$
for all $n \in S^c_{\lambda,\mu}$. Hence, the mapping
\begin{align*}
  E \colon (M_q(\lambda\omega) \otimes_\C M_q(\mu\omega))/N)_{\lambda+\mu-2n} \rarr  (M_q(\lambda\omega) \otimes_\C M_q(\mu\omega))/N)_{\lambda+\mu-2n+2}
\end{align*}
has a nontrivial kernel for \smash{$n \in S^c_{\lambda,\mu}$}, which implies that $(M_q(\lambda\omega) \otimes_\C M_q(\mu\omega))/N$ contains a singular vector $w^{\lambda,\mu}_n$ of weight $q^{(\lambda+\mu-2n)\omega}$ generating a submodule of $(M_q(\lambda\omega) \otimes_\C M_q(\mu\omega))/N$ isomorphic to $M_q((\lambda+\mu-2n)\omega)$ for all $n \in S^c_{\lambda,\mu}$. Since $c_{\lambda+\mu-2n_1} = c_{\lambda+\mu-2n_2}$ for $n_1,n_2 \in S^c_{\lambda,\mu}$ if and only if $n_1=n_2$, we obtain that
\begin{align*}
  (M_q(\lambda\omega) \otimes_\C M_q(\mu\omega))/N \simeq \bigoplus_{n \in S^c_{\lambda,\mu}} M_q((\lambda+\mu-2n)\omega),
\end{align*}
where we used the fact that both sides have the same formal character. Therefore, the short exact sequence \eqref{eq:short exact sequence} can be rewritten into the form
\begin{align*}
  0 \rarr  \bigoplus_{n \in \N_0 \setminus S^c_{\lambda,\mu}} M_q((\lambda+\mu-2n)\omega) \rarr M_q(\lambda\omega) \otimes_\C M_q(\mu\omega) \rarr  \bigoplus_{n \in S^c_{\lambda,\mu}} M_q((\lambda+\mu-2n)\omega) \rarr 0.
\end{align*}
Further, taking the generalized eigenspace decomposition of the $U_q(\mfrak{sl}(2,\C))$-modules in the short exact sequence with respect to the quantum Casimir element $\Cas_q$, we obtain the short exact sequence
\begin{align*}
   0 \rarr  M_q((\lambda+\mu-2n)\omega) \rarr (M_q(\lambda\omega) \otimes_\C M_q(\mu\omega))_{\lambda+\mu-2n}^{\rm cas} \rarr  M_q(-(\lambda+\mu-2n+2)\omega) \rarr 0
\end{align*}
for $n \in S_{\lambda,\mu}$, which is nonsplit since $\dim \Sing (M_q(\lambda\omega) \otimes_\C M_q(\mu\omega))_{\lambda+\mu-2n}^{\rm cas}=2$, and the short exact sequence
\begin{align*}
   0 \rarr  M_q((\lambda+\mu-2n)\omega) \rarr (M_q(\lambda\omega) \otimes_\C M_q(\mu\omega))_{\lambda+\mu-2n}^{\rm cas} \rarr  0 \rarr 0
\end{align*}
for $n \in R_{\lambda,\mu}$, where $(M_q(\lambda\omega) \otimes_\C M_q(\mu\omega))_{\lambda+\mu-2n}^{\rm cas}$ is the generalized eigenspace of $\Cas_q$ with eigenvalue $c_{\lambda+\mu-2n}$ for all $n \in \N_0$.

Since $(M_q(\lambda\omega) \otimes_\C M_q(\mu\omega))_{\lambda+\mu-2n}^{\rm cas}$ belongs to the category $\mcal{O}_q$ for all $n \in \N_0$, we obtain that $(M_q(\lambda\omega) \otimes_\C M_q(\mu\omega))_{\lambda+\mu-2n}^{\rm cas} \simeq P_q((\lambda+\mu-2n)\omega)$ for all $n \in S_{\lambda,\mu}$. Therefore, we have
\begin{align*}
  M_q(\lambda\omega) \otimes_\C M_q(\mu\omega) \simeq \bigoplus_{n \in S_{\lambda,\mu}} P_q((\lambda+\mu-2n)\omega) \oplus \bigoplus_{n \in R_{\lambda,\mu}} M_q((\lambda+\mu-2n)\omega)
\end{align*}
as $U_q(\mfrak{sl}(2,\C))$-modules.

3) If $q^{2\lambda} \in q^{2\N_0}$ and $q^{2\mu} \in q^{2\N_0}$, then the dimension of $\Sing(M_q(\lambda\omega) \otimes_\C M_q(\mu\omega))_{\lambda+\mu-2n}$ is $2$ for $\max\{\lambda_{\rm int}, \mu_{\rm int}\}+1 \leq n \leq \lambda_{\rm int}+\mu_{\rm int}+1$ and $1$ otherwise as follows from Lemma \ref{lem:singular vectors} (ii). Since $c_{\lambda+\mu-2m}= c_{\lambda+\mu-2n}$ for $m>n$ if and only if $m=\lambda_{\rm int}+\mu_{\rm int}-n+1$ for $n=0,1,\dots,\smash{\lfloor {\lambda_{\rm int} + \mu_{\rm int} \over 2} \rfloor}$, we obtain
\begin{align*}
   N = \bigoplus_{n \in T_{\lambda,\mu}} N_{\lambda+\mu-2n} \oplus \bigoplus_{n \in U_{\lambda,\mu}} N_{\lambda+\mu-2n}.
\end{align*}
As $N_{\lambda+\mu-2n}$ is generated by $v^{\lambda,\mu}_{n,\pm}$ for $n \in U_{\lambda,\mu}$ provided $\lambda_{\rm int}+\mu_{\rm int}$ is even and by $v^{\lambda,\mu}_{n,\pm}$ for $n \in U_{\lambda,\mu}$, $n \neq \lfloor {\lambda_{\rm int} + \mu_{\rm int} \over 2} \rfloor +1$ and by $v^{\lambda,\mu}_n$ for $n=\lfloor {\lambda_{\rm int} + \mu_{\rm int} \over 2} \rfloor +1$ provided $\lambda_{\rm int}+\mu_{\rm int}$ is odd, we get $N_{\lambda+\mu-2n} \simeq M_q((\lambda+\mu-2n)\omega)$ for all $n \in U_{\lambda,\mu}$. Further, we have that $N_{\lambda+\mu-2n}$ is generated by $v^{\lambda,\mu}_n$ and \smash{$v^{\lambda,\mu}_{m,+}$}, \smash{$v^{\lambda,\mu}_{m,-}$} with $m=\lambda_{\rm int}+\mu_{\rm int}-n+1$ for $n \in T_{\lambda,\mu} \setminus S_{\lambda,\mu}$. Hence, the submodule of $N_{\lambda+\mu-2n}$ generated by $v^{\lambda,\mu}_n$ is isomorphic to $M_q((\lambda+\mu-2n)\omega)$, which is a reducible $U_q(\mfrak{sl}(2,\C))$-module containing a singular vector \smash{$u^{\lambda,\mu}_{m,1}$} of weight $q^{-(\lambda+\mu-2n+2)\omega}=q^{(\lambda+\mu-2m)\omega}$ with $m=\lambda_{\rm int}+\mu_{\rm int}-n+1$ as follows from Theorem \ref{thm:Verma module decomposition}. Moreover, this singular vector is a linear combination of \smash{$v^{\lambda,\mu}_{m,+}$} and \smash{$v^{\lambda,\mu}_{m,-}$}, since $N_{\lambda+\mu-2n}$ contains all singular vectors of weight $q^{(\lambda+\mu-2m)\omega}$. Furthermore, there exists a linear combination \smash{$u^{\lambda,\mu}_{m,2}$} of \smash{$v^{\lambda,\mu}_{m,+}$} and \smash{$v^{\lambda,\mu}_{m,-}$} such that \smash{$u^{\lambda,\mu}_{m,1}$} and \smash{$u^{\lambda,\mu}_{m,2}$} are linearly independent singular vectors. Therefore, we get $N_{\lambda+\mu-2n} \simeq M_q((\lambda+\mu-2n)\omega) \oplus M_q(-(\lambda+\mu-2n+2)\omega)$ for all $n \in T_{\lambda,\mu} \setminus S_{\lambda,\mu}$. Finally, for $n \in S_{\lambda,\mu}$ we have that $N_{\lambda+\mu-2n}$ is generated by $v^{\lambda,\mu}_n$ and $v^{\lambda,\mu}_m$ with $m=\lambda_{\rm int}+\mu_{\rm int}-n+1$. Hence, the submodule of $N_{\lambda+\mu-2n}$ generated by $v^{\lambda,\mu}_n$ is isomorphic to $M_q((\lambda+\mu-2n)\omega)$, which is a reducible $U_q(\mfrak{sl}(2,\C))$-module containing a singular vector of weight $q^{-(\lambda+\mu-2n+2)\omega}=q^{(\lambda+\mu-2m)\omega}$ with $m=\lambda_{\rm int}+\mu_{\rm int}-n+1$ as follows from Theorem \ref{thm:Verma module decomposition}. Moreover, this singular vector is a multiple of $v^{\lambda,\mu}_m$, since $N_{\lambda+\mu-2n}$ contains all singular vectors of weight $q^{(\lambda+\mu-2m)\omega}$. Therefore, we get $N_{\lambda+\mu-2n} \simeq M_q((\lambda+\mu-2n)\omega)$ for all $n \in S_{\lambda,\mu}$. Further, we have
\begin{align*}
  \ch ((M_q(\lambda\omega) \otimes_\C M_q(\mu\omega))/N) = \sum_{n \in S^c_{\lambda,\mu}} \ch M_q((\lambda+\mu-2n)\omega),
\end{align*}
which gives us $\dim (M_q(\lambda\omega) \otimes_\C M_q(\mu\omega))/N)_{\lambda+\mu-2n}=n-\lfloor {\lambda_{\rm int}+\mu_{\rm int}+1 \over 2} \rfloor$
for all $n \in S^c_{\lambda,\mu}$. Hence, the mapping
\begin{align*}
  E \colon (M_q(\lambda\omega) \otimes_\C M_q(\mu\omega))/N)_{\lambda+\mu-2n} \rarr  (M_q(\lambda\omega) \otimes_\C M_q(\mu\omega))/N)_{\lambda+\mu-2n+2}
\end{align*}
has a nontrivial kernel for \smash{$n \in S^c_{\lambda,\mu}$}, which implies that $(M_q(\lambda\omega) \otimes_\C M_q(\mu\omega))/N$ contains a singular vector $w^{\lambda,\mu}_n$ of weight $q^{(\lambda+\mu-2n)\omega}$ generating a submodule of $(M_q(\lambda\omega) \otimes_\C M_q(\mu\omega))/N$ isomorphic to $M_q((\lambda+\mu-2n)\omega)$ for all $n \in S^c_{\lambda,\mu}$. Since $c_{\lambda+\mu-2n_1} = c_{\lambda+\mu-2n_2}$ for $n_1,n_2 \in S^c_{\lambda,\mu}$ if and only if $n_1=n_2$, we obtain that
\begin{align*}
  (M_q(\lambda\omega) \otimes_\C M_q(\mu\omega))/N \simeq \bigoplus_{n \in S^c_{\lambda,\mu}} M_q((\lambda+\mu-2n)\omega),
\end{align*}
where we used the fact that both sides have the same formal character. Therefore, the short exact sequence \eqref{eq:short exact sequence} can be rewritten into the form
\begin{align*}
  0 \rarr  \bigoplus_{n \in \N_0 \setminus S^c_{\lambda,\mu}} M_q((\lambda+\mu-2n)\omega) \rarr M_q(\lambda\omega) \otimes_\C M_q(\mu\omega) \rarr  \bigoplus_{n \in S^c_{\lambda,\mu}} M_q((\lambda+\mu-2n)\omega) \rarr 0.
\end{align*}
Further, taking the generalized eigenspace decomposition of the $U_q(\mfrak{sl}(2,\C))$-modules in the short exact sequence with respect to the quantum Casimir element $\Cas_q$, we obtain the short exact sequence
\begin{align*}
   0 \rarr  M_q((\lambda+\mu-2n)\omega) \rarr (M_q(\lambda\omega) \otimes_\C M_q(\mu\omega))_{\lambda+\mu-2n}^{\rm cas} \rarr  M_q(-(\lambda+\mu-2n+2)\omega) \rarr 0
\end{align*}
for $n \in S_{\lambda,\mu}$, which is nonsplit since $\dim \Sing((M_q(\lambda\omega) \otimes_\C M_q(\mu\omega))_{\lambda+\mu-2n}^{\rm cas})=2$, the short exact sequence
\begin{align*}
   0 \rarr  M_q((\lambda+\mu-2n)\omega) \rarr (M_q(\lambda\omega) \otimes_\C M_q(\mu\omega))_{\lambda+\mu-2n}^{\rm cas} \rarr  0 \rarr 0
\end{align*}
for $n \in U_{\lambda,\mu}$, and the short exact sequence
\begin{align*}
   0 \rarr  M_q((\lambda+\mu-2n)\omega) \oplus M_q(-(\lambda+\mu-2n+2)\omega) \rarr (M_q(\lambda\omega) \otimes_\C M_q(\mu\omega))_{\lambda+\mu-2n}^{\rm cas} \rarr  0 \rarr 0
\end{align*}
for $n \in T_{\lambda,\mu} \setminus S_{\lambda,\mu}$, where $(M_q(\lambda\omega) \otimes_\C M_q(\mu\omega))_{\lambda+\mu-2n}^{\rm cas}$ is the generalized eigenspace of $\Cas_q$ with eigenvalue $c_{\lambda+\mu-2n}$ for all $n \in \N_0$.

Since $(M_q(\lambda\omega) \otimes_\C M_q(\mu\omega))_{\lambda+\mu-2n}^{\rm cas}$ belongs to the category $\mcal{O}_q$ for all $n \in \N_0$, we obtain that $(M_q(\lambda\omega) \otimes_\C M_q(\mu\omega))_{\lambda+\mu-2n}^{\rm cas} \simeq P_q((\lambda+\mu-2n)\omega)$ for all $n \in S_{\lambda,\mu}$. Therefore, we have
\begin{align*}
  M_q(\lambda\omega) \otimes_\C M_q(\mu\omega) \simeq \bigoplus_{n \in S_{\lambda,\mu}} P_q((\lambda+\mu-2n)\omega) \oplus \bigoplus_{n \in R_{\lambda,\mu}} M_q((\lambda+\mu-2n)\omega)
\end{align*}
as $U_q(\mfrak{sl}(2,\C))$-modules. This completes the proof.}

\
\

\section*{Acknowledgments}
V.\,F.\ is supported in part by CNPq (304467/2017-0) and by Fapesp (2014/09310-5); L.\,K.\ is supported by Capes (88887.137839/2017-00) and J.\,Z.\ is supported by Fapesp (2015/05927-0).



\providecommand{\bysame}{\leavevmode\hbox to3em{\hrulefill}\thinspace}
\providecommand{\MR}{\relax\ifhmode\unskip\space\fi MR }
\providecommand{\MRhref}[2]{%
  \href{http://www.ams.org/mathscinet-getitem?mr=#1}{#2}
}
\providecommand{\href}[2]{#2}

\end{document}